%

%
%
%
\magnification=\magstephalf
\input amstex
\documentstyle{amsppt}
\pagewidth{6.5truein}
\pageheight{8.9truein}
\ifx\refstyle\undefinedRLD\else\refstyle{C}\fi

%

%
%
%
%
%
\count255=\catcode`\!
\catcode`\!=11
\ifx\plot!loaded\relax
   \catcode`\!=\count255 \else\let\plot!loaded=\relax\fi
\chardef\plot!savecc=\count255
%
%
\def\plot!zero{0}
\def\plot!one{1}
\newdimen\plotunitx
\newdimen\plot!unitxu
\def\plot!figscalex{1}
\plotunitx=1truebp
\plot!unitxu=1truebp
\newdimen\plotunity
\newdimen\plot!unityu
\def\plot!figscaley{1}
\plotunity=1truebp
\plot!unityu=1truebp
\def\plot!sepscfac{1}
\newtoks\plot!symbol
\newbox\plot!figurebox
\newdimen\plotcurrx
\newdimen\plotcurry
\newdimen\plotlinewidth
\plotlinewidth=.6bp
\def\plot!savemem{0}
\let\plot!global=\relax
\newif\ifplot
\plottrue
\newif\ifplotPS
\newif\ifplotseparate
\newif\ifplot!infig
\newif\ifplot!trans
\newif\ifplot!loctrans
\newif\ifplot!sepscaled
\newif\ifplot!nowraw
\newif\ifplot!globdef
\newif\ifplot!rawbounds
\newif\ifplot!spseq
\newif\ifplot!stacklock
\newif\ifplot!lhcs
\plot!spseqtrue
%
%
\mathchardef\plot!ci="220E
\mathchardef\plot!bu="220F
\mathchardef\plot!cd="2201
\def\plot!sp{ }
\def\plot!em{}
\let\plot!bg={
\let\plot!eg=}
\def\plot!ni{\prevdepth=-1000pt }
\def\plot!sm#1{{\setbox0=\hbox{#1}\ht0=0pt \dp0=0pt \box0 }}
\def\plot!lo#1\plot!re{\def\plot!bo{#1}\plot!it}
\def\plot!it{\plot!bo \let\plot!ne=\plot!it \else\let\plot!ne=\relax\fi
   \plot!ne}
\let\plot!re=\fi
%
%
\long\def\plot!firstlet#1#2\plot!endofarg{#1}
\let\plot!PStrueold=\plotPStrue
\def\plotPStype#1{
   \plot!PStrueold
   \ifcase#1
      \plotPSfalse
   \or
      \def\plot!local{" }
      \def\plot!global{! }
      \def\plot!rawbegin##1{ps:SDict begin ##1\plot!sp end}
      \def\plot!raw##1{ps:SDict begin ##1\plot!sp end}
      \def\plot!rawend##1{ps:SDict begin ##1\plot!sp end}
      \def\plot!trans{currentpoint /p!s1 2 index def /p!s2 1 index def
         translate exec 0 0 moveto p!s1 neg p!s2 neg translate }
      \def\plot!setorig{}
      \def\plot!rawsetcurr{}
      \def\plot!rawsetorig{currentpoint translate }
      \def\plot!PSfile##1{\includegraphics{##1}}
      \plot!globdeftrue
      \plot!stacklocktrue
      \plot!lhcstrue
   \or
      \def\plot!local{ps:: }
      \def\plot!global{ps::[global] }
      \def\plot!rawbegin##1{ps::[inline,begin] ##1}
      \def\plot!raw##1{ps::[inline] ##1}
      \def\plot!rawend##1{ps::[inline,end] ##1}
      \def\plot!trans{Xpos Ypos translate
         exec Xpos neg Ypos neg translate }
      \def\plot!setorig{Xpos Ypos translate }
      \def\plot!rawsetcurr{Xpos Ypos moveto }
      \def\plot!rawsetorig{Xpos Ypos translate }
      \def\plot!PSfile##1{\plot!rawstart{gsave Xpos Ypos translate}%
         \special{ps: plotfile ##1 inline}\plot!rawfinish{grestore}}
      \plot!rawboundstrue
   \or
      \def\plot!local{empty.ps }
      \def\plot!rawbegin##1{empty.ps ##1}
      \def\plot!raw##1{empty.ps ##1}
      \def\plot!rawend##1{empty.ps ##1}
      \def\plot!trans{pop }
      \def\plot!setorig{}
      \def\plot!rawsetcurr{0 0 moveto }
      \def\plot!rawsetorig{}
      \def\plot!PSfile##1{\special{##1}}
      \plot!spseqfalse
   \else
      \immediate\write16{plotPStype: Illegal type specified -- type set to 0}
      \plotPSfalse
   \fi
   \ifplotPS\ifplot!globdef\else\ifx\plot!global\relax\else
      \ifcase\plot!savemem\or
         \plot!globdeftrue\plot!gwarning
      \or
         \immediate\write16{ }
         \immediate\write16{Do you want me to try to conserve TeX memory?}
         \immediate\write16{(If so, the resulting pages must be printed in %
            order from the beginning.)}
         \message{? }
         \read-16 to\plot!inputstr
         \edef\plot!inputstr{\plot!firstlet\plot!inputstr xx\plot!endofarg}
         \if y\plot!inputstr\plot!globdeftrue\fi
         \if Y\plot!inputstr\plot!globdeftrue\fi
      \fi
   \fi\fi\else\global\let\plot!gwarning=\relax\fi
   \ifplot!globdef\global\let\plot!gwarning=\relax\fi
   \plot!PSinit   
}
\def\plotPSask{
   \immediate\write16{ }
   \immediate\write16{Which DVI-to-PS converter are you using?}
   \immediate\write16{  1. Rokicki dvips (Radical Eye)}
   \immediate\write16{  2. ArborText dvips (DVILASER/PS)}
   \immediate\write16{  3. OzTeX [partial functionality]}
   \immediate\write16{  0. None of the above %
      [PostScript insertions replaced or suppressed]}
   \message{Enter a number: }
   \read-16 to\plot!inputstr
   \edef\plot!inputstr{\plot!firstlet\plot!inputstr xx\plot!endofarg}
   \ifcat0\plot!inputstr\else\def\plot!inputstr{Z}\fi
   \ifnum\expandafter`\plot!inputstr<`0\def\plot!inputstr{99}\fi
   \ifnum\expandafter`\plot!inputstr>`9\def\plot!inputstr{99}\fi
   \plotPStype{\plot!inputstr}
}
\def\plotPStrue{\plotPSask}

\let\plot!septrueold=\plotseparatetrue
\def\plotseparatetrue{\plot!septrueold\plottrue}
\let\plot!sepfalseold=\plotseparatefalse
\def\plotseparatefalse{\plot!sepfalseold\plotnosepscale}
\let\plot!falseold=\plotfalse
\def\plotfalse{\plot!falseold\plotseparatefalse}
\def\plotask{
   \immediate\write16{ }
   \immediate\write16{How do you want figures to appear?}
   \immediate\write16{  1. Normally within text}
   \immediate\write16{  2. On separate pages (blank spaces left in text)}
   \immediate\write16{  0. Not at all (blank spaces left in text)}
   \message{Enter a number: }
   \read-16 to\plot!inputstr
   \if 0\plot!firstlet\plot!inputstr xx\plot!endofarg
      \plotfalse
   \else\if 1\plot!firstlet\plot!inputstr xx\plot!endofarg
      \plottrue \plotseparatefalse
   \else\if 2\plot!firstlet\plot!inputstr xx\plot!endofarg
      \plotseparatetrue
   \else
      \immediate\write16%
         {plotask: Illegal option specified -- figures suppressed}
      \plotfalse
   \fi\fi\fi
}
\gdef\plot!gwarning{
   \immediate\write16{ }
   \immediate\write16{plot WARNING: This DVI-to-PS converter does not %
      fully support global}
   \immediate\write16{ \plot!sp\plot!sp definitions -- printing pages out %
      of order may cause problems.}
   \global\let\plot!gwarning=\relax}
%
%
\def\plot!zerocp{\plotcurrx=0pt\plotcurry=0pt\relax}
{\toksdef\abc=7 \catcode`p=12 \catcode`t=12 \global\abc={pt}}
\expandafter\def\expandafter\plot!unpt\expandafter#\expandafter1\the\toks7{#1 }
\def\plotthedim#1{\expandafter\plot!unpt\the#1 }

\def\plotpdim#1{\plotthedim{#1}\ifnum1000=\mag 1.00375 \else\ifx\plot!unmag
   \plot!em 1003.75 \the\mag\plot!sp div \else 1.00375 \fi\fi div }
\def\plotfdims#1#2{\plotthedim{#1}\plotthedim{\plotunitx}div
   \plotthedim{#2}\plotthedim{\plotunity}div }
\def\plot!recip#1{{\dimen6=1in \dimen8=#1in \dimen0=0sp \global\toks7={}
   \plot!lo
      \count255=\dimen6
      \divide\count255 by\dimen8
      \global\toks7=\expandafter\expandafter\expandafter{\expandafter\the
         \expandafter\toks\expandafter7\expandafter\plot!sp\the\count255}
      \ifdim\dimen0=0sp\global\toks7=\expandafter{\the\toks7.}\fi
      \dimen4=\dimen8
      \multiply\dimen4 by\count255
      \advance\dimen6 by -\dimen4
   \ifdim\dimen0<6sp
      \advance\dimen0 by 1sp
      \multiply\dimen6 by 10
   \plot!re}}
\def\plot!applyssf{\plotunitx=\plot!sepscfac\plot!unitxu
   \plotunity=\plot!sepscfac\plot!unityu}
\def\plot!PSssf{\ifplot!sepscaled\plot!sepscfac\plot!sp dup scale \fi}
\def\plot!errmsg{\errmessage
   {plot error: command allowed only within \string\plotfigbegin
   ...\string\plotfigend}}
\def\plot!rawstart#1{\plot!bg\ifplot!nowraw\special{\plot!raw{#1}}\else
   \special{\plot!rawbegin{#1}}\plot!nowrawtrue\plot!rawboundsfalse\fi}
\def\plot!rawfinish#1{\plot!eg\ifplot!nowraw\special{\plot!raw{#1}}\else
   \special{\plot!rawend{#1}}\fi}
\def\plot!rawstartnull{\ifplot!rawbounds\plot!rawstart{}\else
   \plot!bg\plot!nowrawtrue\fi}
\def\plot!rawfinishnull{\plot!eg\ifplot!rawbounds\ifplot!nowraw
   \else\special{\plot!rawend{}}\fi\fi}
%
%
\def\plotfigscales#1 #2 {\def\plot!figscalex{#1} \plot!unitxu=#1truebp
   \def\plot!figscaley{#2} \plot!unityu=#2truebp \plot!applyssf}
\def\plotfigscale#1 {\plotfigscales{#1} {#1} }
\def\plotrotate#1 {\def\plot!symangle{#1}\plot!loctranstrue\plot!transtrue}
\def\plotscale#1 {\def\plot!symscale{#1}\plot!loctranstrue\plot!transtrue}
\def\plotcancelmag{\def\plot!unmag{\ifnum1000=\mag \else1000 \the\mag\plot!sp
   div dup scale \fi}\plot!loctranstrue\plot!transtrue}
\def\plotresumemag{\def\plot!unmag{}}
\long\def\plotgentrans#1\plotPSend{\plot!loctranstrue\plot!transtrue
   \long\def\plot!transcode{\ifplot!lhcs 1 -1 scale \fi
   #1\ifplot!lhcs\plot!sp 1 -1 scale\fi}}
\def\plotnotrans{\def\plot!transcode{\ifx\plot!symangle\plot!zero\else
   \plot!symangle\plot!sp\ifplot!lhcs neg \fi rotate \fi
   \ifx\plot!symscale\plot!one\else\plot!symscale\plot!sp dup scale \fi}
   \plotresumemag\plot!loctransfalse\ifplot!sepscaled
   \else\plot!transfalse\fi\def\plot!symscale{1}\def\plot!symangle{0}}
\def\plotsepscale#1 {\ifplotseparate\def\plot!sepscfac{#1}
   \ifx\plot!sepscfac\plot!one\plotnosepscale\else
   \plot!sepscaledtrue\plot!transtrue\plot!applyssf \fi\fi}
\def\plotnosepscale{\plot!sepscaledfalse\ifplot!loctrans\else\plot!transfalse
   \fi\def\plot!sepscfac{1}\plot!applyssf}
%
%
\def\plotfigbegin{\setbox\plot!figurebox=\vbox\plot!bg
   \plot!infigtrue\plotnotrans\plot!zerocp}
\def\plotfigend{\xdef\plot!gtemp{\plot!sepscfac}\plot!eg
   {\dimen2=\ht\plot!figurebox
   \dimen4=\wd\plot!figurebox
   \dimen6=\dp\plot!figurebox
   \ifplotseparate
      \shipout\vbox to9truein{\vss\hbox to6.5truein{\hss\box\plot!figurebox
         \hss}\vss}
      \plot!recip{\plot!gtemp}
      \setbox\plot!figurebox=\hbox to \the\toks7\dimen4{\hfil
         \vrule width0pt height\the\toks7\dimen2 depth\the\toks7\dimen6}
   \else\ifplot\else
      \setbox\plot!figurebox=\hbox to \dimen4{\hfil
         \vrule width0pt height\dimen2 depth\dimen6}
   \fi\fi
   \box\plot!figurebox}}
\def\plot!plot{\ifplot
   \vbox to0pt{\kern-\dimen8
   \hbox{%
      \kern\dimen6
      \let\plot!temp=\plot!em
      \ifplotPS\ifplot!trans\ifplot!spseq
         \edef\plot!temp{\plot!unmag\plot!PSssf\plot!transcode}\fi\fi\fi
      \ifx\plot!temp\plot!em\else\plot!rawstart{{\plot!temp} \plot!pla}\fi
      \plot!sm{\hbox to0pt{\hss\the\plot!symbol\hss}}%
      \ifx\plot!temp\plot!em\else\plot!rawfinish{grestore}\fi
   }
   \kern\dimen8}
   \plot!ni\fi}
\def\plot#1 #2 {\ifplot!infig
   {\dimen6=#1\plotunitx\dimen8=#2\plotunity\plot!plot}
   \else\plot!errmsg\fi}

\def\plotmove#1 #2 {\ifplot!infig
   \plotcurrx=#1\plotunitx\plotcurry=#2\plotunity
   \else\plot!errmsg\fi}
\def\plotrmove#1 #2 {\ifplot!infig
   \advance\plotcurrx by #1\plotunitx \advance\plotcurry by #2\plotunity
   \else\plot!errmsg\fi}
%
%
\def\plottext#1{\plot!symbol={\lower\fontdimen22\textfont2\hbox{#1}}}
\def\plotfmla#1{\plottext{$#1$}}
\def\plotsfmla#1{\plot!symbol={\lower\fontdimen22\scriptfont2\hbox
   {$\scriptstyle #1$}}}
\def\plotssfmla#1{\plot!symbol={\lower\fontdimen22\scriptscriptfont2\hbox
   {$\scriptscriptstyle #1$}}}
\def\plotcentered#1{\plot!symbol={\vbox to0pt{\vss\hbox{#1}\vss}}}

%
%
\long\def\plotPSbegin#1\plotPSend{
   \ifplot!infig
   \ifplot\ifplotPS
      \let\plot!temp=\plot!em
      \ifx\plot!one\plot!figscalex\else\let\plot!temp=\plot!one\fi
      \ifx\plot!one\plot!figscaley\else\let\plot!temp=\plot!one\fi
      \ifx\plot!temp\plot!em\else\def\plot!temp{\plot!figscalex\plot!sp
         \ifdim\plot!unitxu=\plot!unityu dup dup scale \else
         \plot!figscaley\plot!sp 2 copy scale dup mul exch
         dup mul 2 div sqrt \fi div }\fi
      \special{\plot!local
         gsave \plot!setorig
         \ifnum1000=\mag\else1000 \the\mag\plot!sp div dup scale \fi
         \ifplot!sepscaled\plot!sepscfac\plot!sp dup scale \fi
         \plotpdim{\plotlinewidth}\plot!temp setlinewidth
         \plotfdims{\plotcurrx}{\plotcurry}moveto #1\plot!sp grestore}
   \fi\fi
   \else\plot!errmsg\fi}
\long\def\plotPSglobal#1\plotPSend{%
   \ifplotPS
      \ifx\plot!global\relax
         \errmessage{plotPSglobal: This PS type does not support global %
            insertions}%
      \else
         \special{\plot!global #1}%
         \plot!gwarning
      \fi
   \fi }
\long\def\plotPSmath#1\plotPSend{{%
   \ifplotPS
      \mathchoice
         {\special{\plot!raw{1.00001}}}%
         {\special{\plot!raw{1.0}}}%
         {\special{\plot!raw{0.7}}}%
         {\special{\plot!raw{0.5}}}
      \special{\plot!local gsave \plot!setorig #1\plot!sp
         \ifplot!stacklock 0 \fi grestore}
      \ifplot!stacklock\special{\plot!raw{count 0 ne {pop} if}}\fi
   \fi}}
\long\def\plotPSgdef#1#2\plotPSend{{%
   \ifplotPS
      \escapechar=-1 %
      \ifplot!globdef
         \special{\plot!global /p!0\string#1 {#2} def}%
         \xdef#1{p!0\string#1 }%
      \else
         \xdef#1{#2\plot!sp}%
      \fi
   \else\global\let#1=\relax\fi}}
\def\plot!befsetup{\setbox0=\hbox{%
   \hskip 1in\raise1in\hbox{\special{\plot!raw{\plot!bfa}}}}%
   \ht0=0pt \wd0=0pt \box0}

\long\def\plotPSbefore#1#2\plotPSend{{%
   \setbox8=\hbox{#1}%
   \ifplotPS
      \hbox{\plot!befsetup\special{\plot!raw{\plot!bfb #2\plot!sp
      grestore}}\box8}%
   \else \box8\fi }}
\long\def\plotPSduring#1#2\plotPSend{{%
   \toks0={{#1}}%
   \ifplotPS\ifplot!spseq
      \hbox{\plot!rawstart{gsave #2}\the\toks0 \plot!rawfinish{grestore}}%
   \else \the\toks0\fi \else \the\toks0\fi }}
\long\def\plotPSbefaft#1#2\plotPSend#3\plotPSend{{%
   \toks0={{#1}}\setbox8=\hbox{\the\toks0}\def\plot!bef{#2}\def\plot!aft{#3}%
   \ifplotPS\ifplot!spseq
      \setbox8=\hbox{\vrule width\wd8 height\ht8 depth\dp8}%
      \hbox to\wd8{\plot!befsetup
         \plot!rawstart{\plot!bab
            \ifx\plot!bef\plot!em\else
            {\plot!bac \plot!bef\plot!sp \plot!bad} \plot!tra \fi}%
         \hbox to0pt{\the\toks0 \hss}%
         \plot!rawfinish{\ifx\plot!aft\plot!em\else {\plot!bac
            \plot!aft} \plot!tra \fi grestore}\hfil}%
   \else \box8\fi \else \box8\fi }}
\long\def\plotPSraw#1\plotPSend{%
   \ifplotPS\ifplot!spseq\special{\plot!raw{#1}}\fi\fi }
%
%
\newdimen\plotdotspc
\plotdotspc=2bp
\def\plotDRline{\ifplot!infig
   \ifplot\ifplotPS
      {\plotfigscale 1
      \edef\plot!litemp{\plotfdims{\dimen0}{\dimen2}}%
      \plotPSbegin\plot!litemp\plot!lia \plotPSend}
   \else\plot!TeXline\fi\fi
   \else\plot!errmsg\fi}
\def\plot!TeXline{{
   \plotlinewidth=\plot!sepscfac\plotlinewidth
   \plotdotspc=\plot!sepscfac\plotdotspc
   \ifdim\dimen0=\plotcurrx
      \dimen6=\dimen0
      \dimen4=\dimen2
      \advance\dimen4 by -\plotcurry
      \ifdim\dimen4<0sp\dimen4=-\dimen4\fi
      \advance\dimen4 by \plotlinewidth
      \plotcentered{\vrule width\plotlinewidth height\dimen4}
      \dimen8=\dimen2
      \advance\dimen8 by \plotcurry
      \divide\dimen8 by 2
      \plot!plot
   \else\ifdim\dimen2=\plotcurry
      \dimen8=\dimen2
      \dimen4=\dimen0
      \advance\dimen4 by -\plotcurrx
      \ifdim\dimen4<0sp\dimen4=-\dimen4\fi
      \advance\dimen4 by \plotlinewidth
      \plotcentered{\vrule width\dimen4 height\plotlinewidth}
      \dimen6=\dimen0
      \advance\dimen6 by \plotcurrx
      \divide\dimen6 by 2
      \plot!plot
   \else\plot!TeXdotline\fi\fi}}
\def\plot!TeXdotline{      
   \plotcentered{\vrule width\plotlinewidth height\plotlinewidth}
   \ifdim\dimen0<\plotcurrx
      \dimen6=\dimen0
      \dimen8=\dimen2
      \dimen0=\plotcurrx
      \dimen2=\plotcurry
   \else
      \dimen6=\plotcurrx
      \dimen8=\plotcurry
   \fi
   \advance\dimen2 by -\dimen8
   \dimen4=\dimen0
   \advance\dimen4 by -\dimen6
   \global\dimen1=\dimen2
   \ifdim\dimen1<0sp\global\dimen1=-\dimen1\fi
   \global\advance\dimen1 by \dimen4
   \global\divide\dimen1 by \plotdotspc
   \global\advance\dimen1 by 1sp
   \global\dimen7=\dimen4
   \global\divide\dimen7 by \dimen1
   \global\dimen3=\dimen7
   \global\multiply\dimen3 by -\dimen1
   \global\advance\dimen3 by \dimen4
   \global\dimen9=\dimen2
   \global\divide\dimen9 by \dimen1
   \global\dimen5=\dimen9
   \global\multiply\dimen5 by -\dimen1
   \global\advance\dimen5 by \dimen2
   \ifdim\dimen5<0sp
      \global\advance\dimen5 by \dimen1
      \global\advance\dimen9 by -1sp
   \fi
   \dimen2=0sp
   \dimen4=0sp
   \plot!lo
      \plot!plot
   \ifdim\dimen6<\dimen0
      \advance\dimen2 by \dimen3
      \advance\dimen4 by \dimen5
      \advance\dimen6 by \dimen7
      \advance\dimen8 by \dimen9
      \ifdim\dimen2<\dimen1\else
         \advance\dimen6 by 1sp
         \advance\dimen2 by -\dimen1
      \fi
      \ifdim\dimen4<\dimen1\else
         \advance\dimen8 by 1sp
         \advance\dimen4 by -\dimen1
      \fi
   \plot!re
}
\def\plotline#1 #2 {
   {\dimen0=#1\plotunitx \dimen2=#2\plotunity
   \plotDRline}
   \plotmove {#1} {#2}
}
\def\plotrline#1 #2 {
   {\dimen0=\plotcurrx \dimen2=\plotcurry
   \advance\dimen0 by #1\plotunitx \advance\dimen2 by #2\plotunity
   \plotDRline}
   \plotrmove {#1} {#2}
}
%
%
\def\plotvskip#1 {\vskip #1\plotunity\relax}
\def\plothphant#1 {\vbox to0pt{\hbox to#1\plotunitx{\hfil}}\plot!ni}
%
%
\def\plot!PSinit{
   \plotPSgdef\plot!circ
      gsave
      \plot!setorig
      newpath
      0 0 2.16 0 360 arc fill
      1 setgray
      newpath
      0 0 1.44 0 360 arc fill
      grestore
   \plotPSend
   \plotPSgdef\plot!bullet
      gsave
      \plot!setorig
      newpath
      0 0 2.16 0 360 arc fill
      grestore
   \plotPSend
   \plotPSgdef\plot!cdot
      gsave
      \plot!setorig
      newpath
      0 0 0.44 0 360 arc fill
      grestore
   \plotPSend
   \plotPSgdef\plot!tra \plot!trans \plotPSend
   \plotPSgdef\plot!pla gsave \plot!tra \plotPSend
   \plotPSgdef\plot!bfa \plot!rawsetcurr currentpoint transform\plotPSend
   \plotPSgdef\plot!bfb gsave \plot!rawsetorig itransform
      72 div exch 72 div exch scale\plotPSend
   \plotPSgdef\plot!bab gsave \plot!rawsetorig itransform
      72 div /p!s4 exch def 72 div /p!s3 exch def grestore gsave\plotPSend
   \plotPSgdef\plot!bac p!s3 p!s4 scale\plotPSend
   \plotPSgdef\plot!bad 1 p!s3 div 1 p!s4 div scale\plotPSend
   \plotPSgdef\plot!lia currentpoint newpath moveto lineto stroke\plotPSend
   \plot!reclaim
}
\def\plot!reclaim{
   \def\plotPSask{\immediate\write16{plot: PStype has already been set}}
   \def\plotPSfalse{\plotPSask}
   \def\plotPStype##1{\plotPSask}
   \let\plot!PSinit=\relax
   \ifplotPS
      \let\plot!TeXline=\relax
      \let\plot!TeXdotline=\relax
   \fi
   \let\plot!reclaim=\relax
}
\catcode`\!=\plot!savecc

\plottrue
\plotPSask
\plotfigscale 7.2

\plotPSgdef\setdoline
   /doline {
      gsave
      12 10 scale
      0.6 72 div setlinewidth
      3 -1 roll 0 translate
      newpath
      2 copy le
      {
         dup 3 1 roll lt
         { .63 .43 }
         { .63 .25 }
         ifelse
         2 index neg moveto
         exch neg lineto
      }
      {
         .15 add
         exch .15 sub
         .25 exch neg moveto
         .65 exch neg lineto
      }
      ifelse
      stroke
      grestore
   } def
\plotPSend

\define\jalg{{{\Cal A}_j}}
\define\jalgp{{{\Cal P}_j}}
\define\crit#1{{cr({#1})}}
\define\threecrit#1/#2/#3/{{\aligned {#1} &{\mapsto} \hbox to0pt{$#2$\hss}
      \hphantom{\alpha_0^{(0)}{=}\alpha_0^{(0)}} \\
      \vspace{-1\jot} &{\mapsto} {#3} \endaligned}}
\define\fourcrit#1/#2/#3/#4/{{\aligned {#1} &{\mapsto} \hbox to0pt{$#2$\hss}
      \hphantom{\alpha_0^{(0)}{=}\alpha_0^{(0)}} \\ \vspace{-1\jot}
      &{\mapsto} {#3} \\ \vspace{-1\jot} &{\mapsto} {#4} \endaligned}}
\define\threeseq#1/#2/#3/#4/{${#1}$~has critical sequence beginning
      ${#2}\mapsto{#3}\mapsto{#4}$}
\define\fourseq#1/#2/#3/#4/#5/{${#1}$~has critical sequence beginning
      ${#2}\mapsto{#3}\mapsto{#4}\mapsto{#5}$}
\define\threeseqv#1/#2/#3/#4/{${#1}$~has a critical vector sequence beginning
      ${#2}\mapsto{#3}\mapsto{#4}$}
\define\fourseqv#1/#2/#3/#4/#5/{${#1}$~has a critical vector sequence beginning
      ${#2}\mapsto{#3}\mapsto{#4}\mapsto{#5}$}
\define\jsub#1{j_{(#1)}}
\define\jsup#1{j^{[#1]}}
\define\et{{\tilde a}}
\define\at{{\tilde \alpha}}
\define\Ct{{\tilde C}}
\define\Ctfunc#1#2{\Ct_{[#1,#2)}}
\define\restrict{\restriction}
\define\Lrestrict{\overset*\to\bigcap}
\define\Lequiv#1/{\overset{#1}\to=}
\define\id{\text{id}}
\define\rightblack{ \null\nobreak\hfill$\blacksquare$}
\define\QED{\rightblack\enddemo}
\define\eps{\varepsilon}

\topmatter
\title Critical points in an algebra of elementary embeddings\endtitle
\author Randall Dougherty\endauthor
\affil Ohio State University\endaffil
\date May 27, 1993 \enddate
\thanks The author was supported by NSF grant number DMS-9158092.\endthanks
\address Department of Mathematics, Ohio State University,
Columbus, OH 43210\endaddress
\abstract
Given two elementary embeddings from the collection of sets of
rank less than $\lambda$ to itself, one can combine them to obtain
another such embedding in two ways: by composition, and by applying one
to (initial segments of) the other.  Hence, a single such nontrivial
embedding $j$ generates an algebra of embeddings via these two operations,
which satisfies certain laws (for example, application distributes over
both composition and application).  Laver has shown, among other things,
that this algebra is free on one generator with respect to these laws.

The set of critical points of members of this algebra is the subject of
this paper.  This set contains the critical point $\kappa_0$ of $j$,
as well as all of the other ordinals $\kappa_n$ in the critical sequence
of $j$ (defined by $\kappa_{n+1} = j(\kappa_n)$).  But the set includes
many other ordinals as well.  The main result of this paper is that
the number of critical points below $\kappa_n$ (which has been shown to be
finite by Laver and Steel) grows so quickly with $n$ that it dominates
any primitive recursive function.  In fact, it grows faster than the
Ackermann function, and even faster than a slow iterate of the Ackermann
function.  Further results show that, even just below $\kappa_4$, one can
find so many critical points that the number is only expressible using
fast-growing hierarchies of iterated functions (six levels of iteration
beyond exponentials).
\endabstract
\endtopmatter

\document

%

%
\head 1. Introduction \endhead

Let $V_\lambda$ be the collection of sets of rank less than $\lambda$, where
$\lambda$ is some fixed limit ordinal.  The assumption that there is a
nontrivial elementary embedding from $V_\lambda$ to $V_\lambda$ is an
extremely strong large cardinal hypothesis~\cite{8}.  But once one such
embedding is known to exist, more of them can be obtained by applying
embeddings to each other.  If $a$ and $b$ are two such embeddings, then
one cannot literally apply $a$ to $b$ since $b \notin V_\lambda$, but
one can apply $a$ to initial segments of $b$, so we define $a(b)$ to
be $\bigcup_{\alpha < \lambda} a(b\restrict V_\alpha)$.  Then $a(b)$
will also be an elementary embedding from $V_\lambda$ to $V_\lambda$.
Of course, one can also obtain new embeddings from old ones by
composition.

Let $\jalg$ be the collection of embeddings generated by a single nontrivial
embedding $j \colon V_\lambda \to V_\lambda$ using the application operation.
Elementarity implies that application is left distributive over itself:
$a(b(c)) = a(b)(a(c))$.  Laver~\cite{5} has shown that $\jalg$ is a free
algebra on the generator $j$ with respect to the left distributive law.
Dehornoy~\cite{2} has carried out an algebraic study of left distributivity;
in particular, Dehornoy~\cite{1} proved a property (irreflexivity) of the free
left distributive algebra on one generator which previously had only been
known under the above large cardinal assumption~\cite{5}.  (A simplification
of the proof of
this particular result has since been given by Larue~\cite{4}.)

If one considers the larger algebra $\jalgp$ generated by $j$ using both
application and composition, then one immediately gets the laws
$a \circ (b \circ c) = (a \circ b) \circ c$ and $(a \circ b)(c) = a(b(c))$,
and elementarity of $a$ gives $a(b \circ c) = a(b) \circ a(c)$ and
$a \circ b = a(b) \circ a$ (the latter because $a(b(x)) = a(b)(a(x))$).
Laver~\cite{5} shows that $\jalgp$ is free with respect to these four laws.

Any embedding $k \in \jalgp$ maps ordinals to ordinals in a strictly increasing
manner, so $k(\alpha) \ge \alpha$ for all $\alpha < \lambda$; the least
$\alpha$ such that $k(\alpha) > \alpha$ is called the {\it critical point} of
$k$, and denoted by $\crit k$.  This $\alpha$ is inaccessible, measurable,
etc., and $k(x) = x$ for $x \in V_\alpha$~\cite{8}.  Elementarity implies that
$k'(\crit k) = \crit{k'(k)}$ and $k'(k(\beta)) = k'(k)(k'(\beta))$; using
these rules, we can obtain many ordinals below $\lambda$ as critical points of
members of $\jalg$.  To start with, let $\kappa_0 = \crit j$ and $\kappa_{n+1}
= j(\kappa_n)$ for $n \in \omega$; the ordinals $\kappa_n$ form a strictly
increasing sequence (called the {\it critical sequence} of $j$), and all of
them are critical points of members of $\jalg$ (if $\kappa_n = \crit k$, then
$\kappa_{n+1} = \crit {j(k)}$).  Perhaps surprisingly, it turns out that other
ordinals also occur as critical points.  To see this, define the
sequence of embeddings $\jsub n$ by $\jsub 1 = j$ and $\jsub {n+1} =
\jsub n(j)$; then we easily compute $\crit{\jsub 2} = j(\kappa_0) =
\kappa_1$, $\jsub 2(\kappa_1) = j(j)(j(\kappa_0)) = j(j(\kappa_0)) =
\kappa_2$, $\jsub 2(\kappa_2) = \kappa_3$, $\crit{\jsub 3} = \jsub
2(\kappa_0) = \kappa_0$, $\jsub 3(\kappa_0) = \jsub 2(j)(\jsub
2(\kappa_0)) = \jsub 2(j(\kappa_0)) = \kappa_2$, and $\jsub 3(\kappa_2)
= \kappa_3$, so $\jsub 3(\kappa_1) = \crit{\jsub 3(\jsub 2)}$ must lie
strictly between $\kappa_2$ and $\kappa_3$.

Further computations yield a number of critical points between
$\kappa_3$ and $\kappa_4$ (for example, the ordinals $j(\jsub
3(\kappa_1))$, $\jsub 2(\jsub 3(\kappa_1))$, and $\jsub 3(\jsub
3(\kappa_1))$ turn out to be distinct).  In fact, one cannot
immediately rule out the possibility that there are infinitely many
critical points below $\kappa_4$.  This turns out not to be the case;
results of Steel and Laver~\cite{6} show that the collection of
critical points of members of $\jalgp$ has order type $\omega$.  (To be
precise, Laver showed that the ordinals $\crit{\jsub n}, n =
1,2,\dotsc$ are the first~$\omega$ critical points of members of
$\jalgp$, listed with duplications, and Steel showed that the ordinals
$\crit{\jsub n}, n = 1,2,\dotsc$ are cofinal in $\lambda$.)

The main result of this paper is that there are very many critical points
other than the ordinals~$\kappa_n$.

\proclaim{Theorem 1} If $F(n)$ is the number of critical points of members of
$\jalgp$ lying below $\kappa_n$, then $F$ grows faster than any primitive
recursive function of $n$. 
\endproclaim

This can be made more precise by considering the standard fast-growing
hierarchy of functions.  (There are slight variations in the ``standard''
definition of this hierarchy.  The version used here is based on a formula of
Hermes and P\'eter, which makes the calculations in this paper simpler; using
another version would not change results like Theorem~1, but it would affect
the explicit formulas relating the numbers of critical points to the
hierarchy.  See Rose~\cite{7} for further details.)  Let $F_0(k) = k+1$, and
define $F_{n+1}$ by iterating $F_n$: $F_{n+1}(0) = F_n(1)$, $F_{n+1}(k+1) =
F_n(F_{n+1}(k))$.  So $F_3(k) = 2^{k+3}-3$, $F_4$ is an iterated exponential,
and so on.  Then diagonalize to obtain $F_\omega(k) = F_k(k)$, and iterate
$F_\omega$ to get $F_{\omega+1}$, etc.  Standard results show that any
primitive recursive function is eventually dominated by $F_n$ for some finite
$n$, and therefore by $F_\omega$ (which is a variant of the Ackermann
function).

An easy induction on terms in $\jalgp$ shows that any element of $\jalgp$
can be written as a composition of one or more members of $\jalg$.  Since
the critical point of a composition of embeddings is the minimum of the
critical points of the individual embeddings, this shows that all critical
points of members of $\jalgp$ are actually critical points of members of
$\jalg$.  (In fact, except in section~3, all of the embeddings used in
the constructions are actually in $\jalg$.)

Section 2 of this paper is a proof that $F(n) > F_\omega(n-1)$ for
all $n \ge 4$; this suffices to prove the above theorem.  In section~3
this result is improved to: $F(n) > F_{\omega+1}(\lfloor \log_3 n \rfloor
- 1)$ for all $n \ge 3$.  These results are strong asymptotically, but do
not say much about $F(n)$ for specific small~$n$.  It turns out that
$F(3) = 4$, but $F(4)$ is quite large; section~4 gives most of a proof
that $F(4)$ exceeds $F_9(F_8(F_8(254)))$, which is an incomprehensibly
large number.  Section~5 covers computation of critical point inequalities,
which leads to some interesting problems; this section can be read
independently of sections~2--4, although one of its purposes is to provide
a few inequalities needed to complete the proof in section~4.

It follows from the theorem above that the algebraic content of the
Steel-Laver result on critical points (see Laver~\cite{6} and Dougherty~and
Jech~\cite{3}) cannot be proved by primitive recursive methods
(i.e., cannot be proved in the standard theory of Primitive Recursive
Arithmetic).  This contrasts with the irreflexivity property of the free
algebra, which, as Dehornoy shows, is provable by primitive recursive
methods.

\head 2. The main construction \endhead

We will now prove Theorem~1 through a sequence of lemmas.  The lemmas will
be self-contained, but the full construction given by them may be hard to
visualize; therefore, after the proof is complete, a picture of the
whole construction and some additional explanatory comments will be given.

Throughout the rest of the paper, `critical point' will mean `critical
point of a member of $\jalgp$' (which, as shown in the preceding
section, is the same as `critical point of a member of $\jalg$'),
and `embedding' will mean `member of $\jalgp$.'

We start by computing the critical sequences of certain simple members
of $\jalg$.  By definition, the critical sequence of $j$ itself is
$\langle \kappa_n \colon n \in \omega\rangle$, or, in a more suggestive
notation, $$\kappa_0 \mapsto \kappa_1 \mapsto \kappa_2 \mapsto \kappa_3
\mapsto \dotsb.$$  Also, if the critical sequence of $i$ is $$\mu_0
\mapsto \mu_1 \mapsto \mu_2 \mapsto \mu_3 \mapsto \dotsb,$$ then the
critical sequence of $k(i)$ is $$k(\mu_0) \mapsto k(\mu_1) \mapsto
k(\mu_2) \mapsto k(\mu_3) \mapsto \dotsb.$$  Therefore, if we define
the sequence of embeddings $\jsup n$ by $\jsup 0 = j$ and $\jsup {n+1}
= j(\jsup n)$, then $\jsup n$~has critical sequence $$\kappa_n \mapsto
\kappa_{n+1} \mapsto \kappa_{n+2} \mapsto \kappa_{n+3} \mapsto
\dotsb.$$

By applying the embeddings $\jsup n$ to each other, we can get
embeddings with slightly more complicated critical sequences.  For
example, the embedding $\jsup1(j)$ has critical sequence $$\kappa_0
\mapsto \kappa_2 \mapsto \kappa_3 \mapsto \kappa_4 \mapsto \dotsb,$$
and the embedding $\jsup3(\jsup3(\jsup1))$ has critical sequence
$$\kappa_1 \mapsto \kappa_2 \mapsto \kappa_5 \mapsto \kappa_6 \mapsto
\kappa_7 \mapsto
\dotsb.$$  In general, given any finite specification
$$\kappa_{n(0)} \mapsto \kappa_{n(1)} \mapsto \dots \mapsto \kappa_{n(l)}$$
with $n(0) < n(1) < \dots < n(l)$, one can find a member of
$\jalg$ whose critical sequence begins with the given
specification: start with $j$, apply $j$ $n(0)$ times to move
the critical point up to $\kappa_{n(0)}$, then apply $\jsup{n(0)+1}$
$n(1)-n(0)-1$ times to move the next member of the critical sequence
up to $\kappa_{n(1)}$, then apply $\jsup{n(1)+1}$ $n(2)-n(1)-1$ times,
and so on.  The examples below will only use some simple cases
of this; for instance, if one wants an embedding $k$ with critical
sequence beginning $\kappa_2 \mapsto \kappa_n \mapsto \kappa_{n+1}$
($n \ge 3$), then one can let $k = (\jsup 3)^{n-3}(\jsup 2)$
(i.e., start with $\jsup 2$ and apply $\jsup 3$ $n-3$ times).

A main fact to be used in the construction is that each embedding in
$\jalgp$ is a strictly order-preserving map from critical points to
critical points; this will be used over and over again to produce more
and more critical points, and we will similarly get more and more
embeddings as well.  To start with:

\proclaim{Lemma 2} Suppose that one has critical points
$\zeta < \alpha_0 < \alpha_1 < \dots < \alpha_n$ and embeddings
$a_i \in \jalgp$ for $i < n$ such that \threeseq a_i/\zeta/\alpha_i
/\alpha_{i+1}/ for each $i$.  Then there are at least $2^n$ critical
points in the half-open interval\/ $[\zeta,\alpha_n)$. \endproclaim

\demo{Proof} Induct on $n$; for $n=0$, the single critical point
$\zeta$ will suffice.  Given the result for $n$, if one has
$\alpha_i$ for $i \le n+1$ and $a_i$ for $i < n+1$ as above,
then the induction hypothesis gives at least $2^n$ critical
points in $[\zeta,\alpha_n)$.  The embedding $a_n$ maps
the interval $[\zeta,\alpha_n)$ to the interval
$[\alpha_n,\alpha_{n+1})$, and maps distinct critical points in
the former interval to distinct critical points in the latter.  Therefore,
there are at least $2^n$ critical points in $[\alpha_n,\alpha_{n+1})$,
giving a total of at least $2^{n+1}$ critical points in $[\zeta,
\alpha_{n+1})$.  This completes the induction. \QED

For example, if we let $\zeta=\kappa_0$, $\alpha_i = \kappa_{i+1}$, and
$a_i = (\jsup1)^i(j)$, then we get at least $2^n$ critical points below
$\kappa_{n+1}$ for any $n$.  However, we can do much better than this.

\proclaim{Lemma 3} Suppose that one has critical points
$\zeta < \zeta' < \beta_0 < \beta_1 < \dots < \beta_n$ and embeddings
$b_i \in \jalgp$ for $i < n$ such that \threeseq b_i/\zeta'/\beta_i
/\beta_{i+1}/ for each $i$.  Also suppose that there is an embedding
$z$ such that \threeseq z/\zeta/\zeta'
/\beta_0/.  Then there exist critical points $\alpha_i$ ($i \le 2^n$)
and embeddings $a_i$ ($i < 2^n$) such that \threeseq a_i/\zeta/\alpha_i
/\alpha_{i+1}/ for each $i$, $\alpha_0 = \zeta'$, and $\alpha_{2^n} =
\beta_n$.\endproclaim

\demo{Proof} Induct on $n$; for $n=0$, just let $a_0 = z$, $\alpha_0 =
\zeta'$, and $\alpha_1 = \beta_0$.  Given the result for $n$, if one has
$\beta_i$ for $i \le n+1$ and $b_i$ for $i < n+1$ as above,
then the induction hypothesis gives critical points $\alpha_i$ ($i \le 2^n$)
and embeddings $a_i$ ($i < 2^n$) as required, with $\alpha_0 = \zeta'$
and $\alpha_{2^n} = \beta_n$.  Define $\alpha_{2^n+i}$ for $1 \le i \le 2^n$
by the equation $\alpha_{2^n+i} = b_n(\alpha_i)$; note that this
equation also holds for $i = 0$, since $b_n(\zeta') = \beta_n$.
Also note that $b_n(\zeta) = \zeta$, since $\zeta < \zeta' =
\crit{b_n}$.
Hence, if we let $a_{2^n+i} = b_n(a_i)$ for $0 \le i < 2^n$, then
\threeseq a_{2^n+i}/\zeta/\alpha_{2^n+i}/\alpha_{2^n+i+1}/ for
each such $i$.  Finally, $\alpha_{2^{n+1}} = b_n(\beta_n) = \beta_{n+1}$.
This completes the induction. \QED

Note that it follows from the fact that \threeseq a_i/\zeta/\alpha_i
/\alpha_{i+1}/ that $\zeta < \alpha_i < \alpha_{i+1}$.

One could apply this lemma with $\zeta = \kappa_0$, $\zeta' = \kappa_1$,
$\beta_i = \kappa_{i+2}$, $z = j$, and $b_i = (\jsup2)^i(\jsup1)$ to
get $2^n$ embeddings as in Lemma~2 ending up with $\alpha_{2^n} =
\kappa_{n+2}$; hence, Lemma~2 would give $2^{2^n}$ critical points
below $\kappa_{n+2}$.  In fact, for any $m$, if one started with
embeddings $b_i$ ($i \le n$) such that \threeseq
b_i/\kappa_m/\kappa_{m+1+i}/\kappa_{m+2+i}/, then one could apply
Lemma~3 $m$ times and then apply Lemma~2 to get $\exp_2^{m+1}(n)$ (a
tower of exponentials consisting of $m+1$ $2$'s with an $n$ at the
top) critical points below $\kappa_{n+m+1}$.  Again, however, this
is not as much as following results will give.

The next result uses a numerical function $C_2(n)$ defined recursively
as follows: $C_2(0) = 0$ and $C_2(n+1) = C_2(n)+2^{C_2(n)}$.  (So
$C_2$ is something like an iterated exponential.)

\proclaim{Lemma 4} Suppose that one has critical points $\zeta <
\zeta' < \delta_0 < \delta_1 < \dots < \delta_n$ and embeddings $d_i
\in \jalgp$ for $i < n$ such that \fourseq
d_i/\zeta/\zeta'/\delta_i /\delta_{i+1}/ for each $i$.  Then there
exist critical points $\beta_i$ ($i \le C_2(n)$) and embeddings $b_i$ ($i
< C_2(n)$) such that \threeseq b_i/\zeta'/\beta_i /\beta_{i+1}/ for
each $i$, $\beta_0 = \delta_0$, and $\beta_{C_2(n)} =
\delta_n$.\endproclaim

\demo{Proof} Induct on $n$; the case $n=0$ is essentially vacuous (just
let $\beta_0 = \delta_0$).  Given the result for $n$, if one has
$\delta_i$ for $i \le n+1$ and $d_i$ for $i < n+1$ as above, then the
induction hypothesis gives critical points $\beta_i$ ($i \le C_2(n)$)
and embeddings $b_i$ ($i < C_2(n)$) as required, with $\beta_0 =
\delta_0$ and $\beta_{C_2(n)} = \delta_n$.  We can now apply Lemma~3,
using $z = d_0$, to get critical points $\alpha_i$ ($i \le m$) and
embeddings $a_i$ ($i < m$), where $m = 2^{C_2(n)}$, such that \threeseq
a_i/\zeta/\alpha_i /\alpha_{i+1}/ for each $i$, $\alpha_0 =
\zeta'$, and $\alpha_{m} = \beta_{C_2(n)} = \delta_n$.  Define
$\beta_{C_2(n)+i}$ for $1 \le i \le m$ by the equation
$\beta_{C_2(n)+i} = d_n(\alpha_i)$; note that this equation also holds
for $i = 0$, since $d_n(\zeta') = \delta_n$.  Hence, if we let
$b_{C_2(n)+i} = d_n(a_i)$ for $0 \le i < m$, then \threeseq
b_{C_2(n)+i}/\zeta'/\beta_{C_2(n)+i}/\beta_{C_2(n)+i+1}/ for each
such $i$.  Finally, $\beta_{C_2(n+1)} = d_n(\delta_n) = \delta_{n+1}$.
This completes the induction. \QED

\proclaim{Lemma 5} Suppose that one has critical points $\zeta <
\zeta' < \zeta'' < \eps_0 < \eps_1 < \dots < \eps_n$ and embeddings
$e_i \in \jalgp$ for $i < n$ such that \threeseq e_i/\zeta''/\eps_i
/\eps_{i+1}/ for each $i$.  Also suppose that there is an embedding $z$
such that \fourseq z/\zeta/\zeta'/\zeta'' /\eps_0/.  Then there
exist critical points $\delta_i$ ($i \le 2^n$) and embeddings $d_i$ ($i
< 2^n$) such that \fourseq d_i/\zeta/\zeta'/\delta_i /\delta_{i+1}/ for
each $i$, $\delta_0 = \zeta''$, and $\delta_{2^n} =
\eps_n$.\endproclaim

\demo{Proof} Just as in the proof of Lemma 3, let $d_0 = z$ and
$d_{2^n+i} = e_n(d_i)$ for $0 \le i < 2^n$. \QED

We need one more lemma of this sort, for which we need a family
of numerical functions to be denoted $C_{2N+2}(n)$ for $N \ge 1$.
These are again defined recursively, both on $N$ and $n$, by
the formulas $C_{2N+2}(0) = 0$ and $C_{2N+2}(n+1) = C_{2N+2}(n) +
C_{2N}(2^{C_{2N+2}(n)})$.  (So $C_{2N+2}$ is an augmented iteration of
$C_{2N}$.)

\proclaim{Lemma 6} Let $N \ge 1$, and suppose that one has critical
points $\zeta_0 < \zeta_1 < \dots < \zeta_{N+1} < \eta_0 < \eta_1
< \dots < \eta_n$ and embeddings $h_i \in \jalgp$ for $i < n$ such that
\fourseq h_i/\zeta_N/\zeta_{N+1}/\eta_i /\eta_{i+1}/ for each $i$.
Also suppose that there exist embeddings $z_k$ ($k < N$) such that
\fourseq z_k/\zeta_k/\zeta_{k+1}/\zeta_{k+2} /\zeta_{k+3}/ for
each $k < N$, where $\zeta_{N+2} = \eta_0$.  Then there exist
critical points $\eps_i$ ($i \le C_{2N+2}(n)$) and embeddings $e_i$ ($i
< C_{2N+2}(n)$) such that \threeseq e_i/\zeta_{N+1}/\eps_i
/\eps_{i+1}/ for each $i$, $\eps_0 = \eta_0$, and $\eps_{C_{2N+2}(n)} =
\eta_n$.\endproclaim

\demo{Proof} Proceed by induction on $N$; we will prove the result for
$N$, assuming it is true for $N'<N$.  The induction step for $N$ is
similar to the proof of Lemma~4.
 Induct on $n$; the case $n=0$ is essentially vacuous (just let $\eps_0
= \eta_0$).  Given the result for $n$, if one has $\eta_i$ for $i \le
n+1$ and $h_i$ for $i < n+1$ as above, then the induction hypothesis
gives critical points $\eps_i$ ($i \le C_{2N+2}(n)$) and embeddings
$e_i$ ($i < C_{2N+2}(n)$) as required, with $\eps_0 = \eta_0$ and
$\eps_{C_{2N+2}(n)} = \eta_n$.  We can now apply Lemma~5, with $\zeta
= \zeta_{N-1}$, $\zeta'=\zeta_N$, $\zeta'' = \zeta_{N+1}$,
and $z = z_{N-1}$, to get critical points $\delta_i$ ($i \le m$) and
embeddings $d_i$ ($i < m$), where $m = 2^{C_{2N+2}(n)}$, such that
\fourseq d_i/\zeta_{N-1}/\zeta_N/\delta_i /\delta_{i+1}/ for each
$i$, $\delta_0 = \zeta_{N+1}$, and $\delta_{m} = \eps_{C_{2N+2}(n)} =
\eta_n$.  Finally, apply either Lemma~4 (if $N=1$) or the induction
hypothesis for $N-1$ (if $N > 1$) to get critical points $\beta_i$ ($i
\le C_{2N}(m)$) and embeddings $b_i$ ($i < C_{2N}(m)$) such that
\threeseq b_i/\zeta_N/\beta_i /\beta_{i+1}/ for each $i$, $\beta_0 =
\delta_0 = \zeta_{N+1}$, and $\beta_{C_{2N}(m)} = \delta_m =
\eta_n$.  Define $\eps_{C_{2N+2}(n)+i}$ for $1 \le i \le C_{2N}(m)$ by
the equation $\eps_{C_{2N+2}(n)+i} = h_n(\beta_i)$; note that this
equation also holds for $i = 0$, since $h_n(\zeta_{N+1}) = \eta_n$.
Hence, if we let $e_{C_{2N+2}(n)+i} = h_n(b_i)$ for $0 \le i <
C_{2N}(m)$, then \threeseq e_{C_{2N+2}(n)+i}/ \zeta_{N+1}/
\eps_{C_{2N+2}(n)+i}/ \eps_{C_{2N+2}(n)+i+1}/ for each such $i$.
Also, $\eps_{C_{2N+2}(n+1)} = h_n(\eta_n) = \eta_{n+1}$.  This
completes both inductions. \QED

We are now ready to put everything together to get a lower bound
on $F(n)$, the number of critical points below $\kappa_n$.  Let us
define $C_k(n)$ to be $2^n$ if $n$ is not even and positive.  In other
words, the full definition of $C_k(n)$ is:
$$\alignat 2
&\text{If $k$ is $0$ or odd:} & \qquad C_k(n) &= 2^n \\
&\text{If $k=2$:} & \qquad C_k(0) &= 0 \\
& & \qquad C_k(n+1) &= C_k(n) +
      2^{C_k(n)}\\
&\text{If $k$ is even and greater than $2$:} & \qquad C_k(0) &= 0 \\
& & \qquad C_k(n+1) &= C_k(n) +
      C_{k-2}(2^{C_k(n)})
\endalignat$$

\proclaim{Proposition 7} Let $N \ge 0$, and suppose that one has critical
points $\zeta_0 < \zeta_1 < \dots < \zeta_{N+1} < \eta_0 < \eta_1
< \dots < \eta_n$ and embeddings $h_i \in \jalgp$ for $i < n$ such that
\fourseq h_i/\zeta_N/\zeta_{N+1}/\eta_i /\eta_{i+1}/ for each $i$.
Also suppose that there exist embeddings $z_k$ ($k < N$) such that
\fourseq z_k/\zeta_k/\zeta_{k+1}/\zeta_{k+2} /\zeta_{k+3}/ for
each $k < N$, where $\zeta_{N+2} = \eta_0$.  Then the number of
critical points below $\eta_n$ is at least
$C_0(C_1(\dots(C_{2N+2}(n))\dots))$.\endproclaim

\demo{Proof} Induct on $N$.  For $N = 0$, the case $n = 0$ is trivial:
$C_0(C_1(C_2(n)) = 2$ and $\zeta_0 < \zeta_1 < \eta_0$.  If $N=0$
and $n>0$, apply Lemma~4 to get embeddings $b_i$ and critical points
$\beta_i$ so that \threeseq b_i/ \zeta_1/ \beta_i/ \beta_{i+1}/ and
$\beta_0 = \eta_0$, $\beta_{C_2(n)} = \eta_n$; then apply Lemma~3 (with
$z = h_0$) to get embeddings $a_i$ and critical points $\alpha_i$ so
that \threeseq a_i/ \zeta_0/ \alpha_i/ \alpha_{i+1}/ and
$\alpha_{C_1(C_2(n))} = \eta_n$; then apply Lemma~2 to get
$C_0(C_1(C_2))$ critical points below $\eta_n$.  For the induction step
$N>0$, given $\eta_i$ and $h_i$ as above, apply Lemma~6 to get
embeddings $e_i$ and critical points $\eps_i$ so that \threeseq e_i/
\zeta_{N+1}/ \eps_i/ \eps_{i+1}/ and $\eps_0 = \eta_0$,
$\eps_{C_{2N+2}(n)} = \eta_n$; then apply Lemma~5 (with $z = z_{N-1}$)
to get embeddings $d_i$ and critical points $\delta_i$ so that \fourseq
d_i/ \zeta_{N-1}/ \zeta_N/ \delta_i/ \delta_{i+1}/ and $\delta_0 =
\zeta_{N+1}$, $\delta_{C_{2N+1}(C_{2N+2}(n))} = \eta_n$; finally,
apply the induction hypothesis for $N-1$ to get
$C_0(C_1(\dots(C_{2N+2}(n))\dots))$ critical points below~$\eta_n$.
\QED

Now define the function $g$ as follows: if $m \le 2$,
let $g(m) = m$; if $m \ge 3$, let $g(m) =
C_0(C_1(\dots(C_{2m-4}(1))\dots))$.

\proclaim{Corollary 8} Suppose one has critical points
$\zeta_0 < \zeta_1 < \dots < \zeta_m$ and embeddings
$z_k$ ($0 \le k \le m-3$) such that \fourseq z_k/ \zeta_k/
\zeta_{k+1}/ \zeta_{k+2}/ \zeta_{k+3}/.  Then the number
of critical points below $\zeta_m$ is at least $g(m)$. \endproclaim

\demo{Proof} For $m < 3$, this is trivial.  For $m \ge 3$, apply
Proposition~7 with $N = m-3$, $n = 1$, $\eta_0 = \zeta_{m-1}$,
$\eta_1 = \zeta_m$, and $h_0 = z_{m-3}$. \QED

Clearly we can apply Corollary 8 with $\zeta_k = \kappa_k$ and $z_k =
\jsup k$ to get $F(m) \ge g(m)$ for all $m$.  To complete the proof of
Theorem~1, it will suffice to show that $g(m) \ge F_\omega(m-1)$ for $m
\ge 5$.  (In fact, this is true for $m = 4$ as well: direct computation
shows that $g(4) = 256$ and $F_\omega(3) = 61$.)

It is easy to prove by induction on $N$ that $C_{2N}(n)$ is
a strictly increasing function of $n$ for each fixed $N$, and hence
$C_{2N}(n) \ge n$.  Of course, the same statements trivially
hold for $C_{2N+1}$.

We can now prove by double induction on $N$ and $n$ that
$C_{2N}(n+4) > F_{N+3}(n)+3$ for all $N>0$ and $n$.
For $N=1$ and $n=0$, compute directly that $C_2(4) = 2059 > 16 = F_4(0)+3$.
If $N=1$ and $n > 0$, then we have $$C_2(n+4) = C_2(n+3) +2^{C_2(n+3)}
> 2^{C_2(n+3)} > 2^{F_4(n-1) + 3} = F_4(n)+3.$$
For $N>1$ and $n=0$, we have
$$C_{2N}(4) = C_{2N}(3) + C_{2N-2}(2^{C_{2N}(3)})
> C_{2N-2}(2^3) > C_{2N-2}(5) > F_{N+2}(1)+3
= F_{N+3}(0)+3.$$  Finally, if $N>1$ and $n>0$, then
$$\multline C_{2N}(n+4) = C_{2N}(n+3) + C_{2N-2}(2^{C_{2N}(n+3)})
> C_{2N-2}(C_{2N}(n+3)+1) \\
> C_{2N-2}(F_{N+3}(n-1)+4) > F_{N+2}(F_{N+3}(n-1))+3 = F_{N+3}(n)+3.
\endmultline$$

Next, one can show by easy inductions on $N$ that
$C_{2N}(1) = 1$ and $C_{2N}(2) = N+2$ for all $N \ge 1$.  We are now
ready to prove the desired bound on $g(m)$ for $m \ge 5$:
$$\align g(m) &= C_0(C_1(\dots(C_{2m-4}(1))\dots)) \\
&\ge C_{2m-8}(C_{2m-7}(C_{2m-6}(C_{2m-5}(C_{2m-4}(1))))) \\
&= C_{2m-8}(2^{m-1}) \\
&> C_{2m-8}(m+3) > F_{m-1}(m-1) = F_\omega(m-1).\endalign\nopagebreak$$
This completes the proof of Theorem~1.

\medpagebreak

The following remarks may help to clarify what is going on in
the proof of Theorem~1.  (Or they may make it more confusing.  Read
on at your own risk.)

First, let us observe a property which is not hard to see from the proofs of
Lemmas~3--6 but which does not come out in the statements of these
lemmas.  Lemma~3 states that, given a sequence of embeddings $b_i$ ($i <
n$) such that \threeseq b_i/\zeta'/\beta_i/\beta_{i+1}/ (and an
additional embedding $z$), one can obtain embeddings $a_i$ ($i < 2^n$)
such that \threeseq a_i/\zeta/\alpha_i/\alpha_{i+1}/, where
$\alpha_{2^n} = \beta_n$.  Now, if one can obtain some additional
embeddings $b_n,b_{n+1},\dotsc$ to extend the given sequence of $b_i$'s,
then one gets a longer sequence of embeddings $a_i$, and it is
clear from the proof of Lemma~3 that this longer
sequence will extend the sequence of $2^n$ $a_i$'s obtained from
$b_0,\dots,b_{n-1}$.  Applying this repeatedly, we see that, if $2^m$
embeddings $a_i$ are obtained from $m$ embeddings $b_i$, then we will
actually have $\alpha_{2^i} = \beta_i$ for all $i \le m$, not just for
$i = m$.  Also, if we somehow obtained an infinite sequence of embeddings
$b_i$, then we would get an infinite sequence of embeddings $a_i$ and we
would have $\alpha_{2^i} = \beta_i$ for all $i$.  Similar statements
hold for the situations in Lemmas~4--6.

Let us now introduce a more uniform notation.  Instead of referring to
embeddings $a_i,b_i$ and critical points $\zeta,\zeta',\alpha_i,
\beta_i$ in Lemma~3, one can instead state the result as follows:
given embeddings $a_i^{(2)}$ ($i < n$) such that \threeseq a_i^{(2)}/
\zeta_1/\alpha_i^{(2)}/\alpha_{i+1}^{(2)}/, and given that \threeseq
z_0/\zeta_0/\zeta_1/\alpha_0^{(2)}/, one can produce embeddings
$a_i^{(1)}$ for $i < 2^n$ such that \threeseq a_i^{(1)}/
\zeta_1/\alpha_i^{(1)}/\alpha_{i+1}^{(1)}/, where $\alpha_0^{(1)} =
\zeta_1$ and $\alpha_{2^n}^{(1)} = \alpha_n^{(2)}$ (in fact, by the
preceding paragraph, $\alpha_{2^i}^{(1)} = \alpha_i^{(2)}$ for all $i
\le n$).  Similarly, one can write Lemma~4 so that it produces
embeddings $a_i^{(2)}$ and critical points $\alpha_i^{(2)}$ from given
embeddings $a_i^{(3)}$ and critical points~$\alpha_i^{(3)}$ (writing
the two fixed ordinals as $\zeta_0,\zeta_1$ instead of
$\zeta,\zeta'$).  Lemma~5 can be stated so that one uses
$a_i^{(2N+2)}$ and
$\alpha_i^{(2N+2)}$ to produce $a_i^{(2N+1)}$ and $\alpha_i^{(2N+1)}$,
given fixed ordinals $\zeta_{N-1},\zeta_N, \zeta_{N+1}$ and an
initial embedding $z_{N-1}$.  (One needs a variable superscript because
Lemma~5 is applied many times, once for each inductive step in the
proofs of Lemma~6 and Proposition~7.) Finally, Lemma~6 produces
embeddings $a_i^{(2N+2)}$ and critical points $\alpha_i^{(2N+2)}$ from
given embeddings $a_i^{(2N+3)}$ and critical points
$\alpha_i^{(2N+3)}$.

\pageinsert
\vfill
\plotfigbegin
   \normalbaselines
   \vskip 6.7truein
   \plotPSbegin
      \setdoline
      -9 63 translate
      1 1 0 doline
      1 2 1 doline
      1 3 1 doline
      1 4 2 doline
      1 5 2 doline
      1 6 2 doline
      1 7 2 doline
      2 -1 0 doline
      2 1 1 doline
      2 2 1 doline
      2 3 2 doline
      2 4 2 doline
      2 5 2 doline
      2 6 2 doline
      2 7 2 doline
      3 1 0 doline
      3 2 1 doline
      3 3 1 doline
      3 4 2 doline
      3 5 2 doline
      3 6 2 doline
      3 7 2 doline
      4 -1 0 doline
      4 1 1 doline
      4 2 1 doline
      4 3 1 doline
      4 4 2 doline
      4 5 2 doline
      4 6 2 doline
      4 7 2 doline
      5 1 0 doline
      5 2 1 doline
      5 3 1 doline
      5 4 2 doline
      5 5 2 doline
      5 6 2 doline
      5 7 2 doline
   \plotPSend
   \plotfmla{\boxed{\threecrit\zeta_0/\alpha_0^{(1)}{=}\zeta_1/
         \alpha_1^{(1)}{=}\alpha_0^{(2)}/}}
   \plot 3 63
   \plotfmla{\threecrit\zeta_0/\alpha_1^{(1)}/
         \alpha_2^{(1)}{=}\alpha_1^{(2)}/}
   \plot 3 53
   \plotfmla{\threecrit\zeta_0/\alpha_2^{(1)}/\alpha_3^{(1)}/}
   \plot 3 43
   \plotfmla{\threecrit\zeta_0/\alpha_3^{(1)}/
         \alpha_4^{(1)}{=}\alpha_2^{(2)}/}
   \plot 3 33
   \plotfmla{\threecrit\zeta_0/\alpha_4^{(1)}/\alpha_5^{(1)}/}
   \plot 3 23
   \plotfmla{\threecrit\zeta_0/\alpha_5^{(1)}/\alpha_6^{(1)}/}
   \plot 3 13
   \plotfmla{\threecrit\zeta_0/\alpha_6^{(1)}/\alpha_7^{(1)}/}
   \plot 3 3
   \plotfmla{\threecrit\zeta_0/\alpha_7^{(1)}/
         \alpha_8^{(1)}{=}\alpha_3^{(2)}/}
   \plot 3 -7
   \plotfmla{\threecrit\zeta_1/\alpha_0^{(2)}{=}\alpha_0^{(3)}/
         \alpha_1^{(2)}{=}\alpha_1^{(3)}/}
   \plot 15 63
   \plotfmla{\threecrit\zeta_1/\alpha_1^{(2)}/\alpha_2^{(2)}/}
   \plot 15 53
   \plotfmla{\threecrit\zeta_1/\alpha_2^{(2)}/
         \alpha_3^{(2)}{=}\alpha_2^{(3)}/}
   \plot 15 43
   \plotfmla{\threecrit\zeta_1/\alpha_3^{(2)}/\alpha_4^{(2)}/}
   \plot 15 33
   \plotfmla{\threecrit\zeta_1/\alpha_4^{(2)}/\alpha_5^{(2)}/}
   \plot 15 23
   \plotfmla{\threecrit\zeta_1/\alpha_5^{(2)}/\alpha_6^{(2)}/}
   \plot 15 13
   \plotfmla{\threecrit\zeta_1/\alpha_6^{(2)}/\alpha_7^{(2)}/}
   \plot 15 3
   \plotfmla{\threecrit\zeta_1/\alpha_7^{(2)}/\alpha_8^{(2)}/}
   \plot 15 -7
   \plotfmla{\boxed{\fourcrit\zeta_0/\zeta_1/\alpha_0^{(3)}{=}\zeta_2/
         \alpha_1^{(3)}{=}\alpha_0^{(4)}/}}
   \plot 27 63
   \plotfmla{\fourcrit\zeta_0/\zeta_1/\alpha_1^{(3)}/
         \alpha_2^{(3)}{=}\alpha_1^{(4)}/}
   \plot 27 53
   \plotfmla{\fourcrit\zeta_0/\zeta_1/\alpha_2^{(3)}/\alpha_3^{(3)}/}
   \plot 27 43
   \plotfmla{\fourcrit\zeta_0/\zeta_1/\alpha_3^{(3)}/
         \alpha_4^{(3)}{=}\alpha_2^{(4)}/}
   \plot 27 33
   \plotfmla{\fourcrit\zeta_0/\zeta_1/\alpha_4^{(3)}/\alpha_5^{(3)}/}
   \plot 27 23
   \plotfmla{\fourcrit\zeta_0/\zeta_1/\alpha_5^{(3)}/\alpha_6^{(3)}/}
   \plot 27 13
   \plotfmla{\fourcrit\zeta_0/\zeta_1/\alpha_6^{(3)}/\alpha_7^{(3)}/}
   \plot 27 3
   \plotfmla{\fourcrit\zeta_0/\zeta_1/\alpha_7^{(3)}/
         \alpha_8^{(3)}{=}\alpha_3^{(4)}/}
   \plot 27 -7
   \plotfmla{\threecrit\zeta_2/\alpha_0^{(4)}{=}\alpha_0^{(5)}/
         \alpha_1^{(4)}{=}\alpha_1^{(5)}/}
   \plot 39 63
   \plotfmla{\threecrit\zeta_2/\alpha_1^{(4)}/\alpha_2^{(4)}/}
   \plot 39 53
   \plotfmla{\threecrit\zeta_2/\alpha_2^{(4)}/\alpha_3^{(4)}/}
   \plot 39 43
   \plotfmla{\threecrit\zeta_2/\alpha_3^{(4)}/
         \alpha_4^{(4)}{=}\alpha_2^{(5)}/}
   \plot 39 33
   \plotfmla{\threecrit\zeta_2/\alpha_4^{(4)}/\alpha_5^{(4)}/}
   \plot 39 23
   \plotfmla{\threecrit\zeta_2/\alpha_5^{(4)}/\alpha_6^{(4)}/}
   \plot 39 13
   \plotfmla{\threecrit\zeta_2/\alpha_6^{(4)}/\alpha_7^{(4)}/}
   \plot 39 3
   \plotfmla{\threecrit\zeta_2/\alpha_7^{(4)}/\alpha_8^{(4)}/}
   \plot 39 -7
   \plotfmla{\boxed{\fourcrit\zeta_1/\zeta_2/\alpha_0^{(5)}{=}\zeta_3/
         \alpha_1^{(5)}{=}\alpha_0^{(6)}/}}
   \plot 51 63
   \plotfmla{\fourcrit\zeta_1/\zeta_2/\alpha_1^{(5)}/
         \alpha_2^{(5)}{=}\alpha_1^{(6)}/}
   \plot 51 53
   \plotfmla{\fourcrit\zeta_1/\zeta_2/\alpha_2^{(5)}/\alpha_3^{(5)}/}
   \plot 51 43
   \plotfmla{\fourcrit\zeta_1/\zeta_2/\alpha_3^{(5)}/
         \alpha_4^{(5)}{=}\alpha_2^{(6)}/}
   \plot 51 33
   \plotfmla{\fourcrit\zeta_1/\zeta_2/\alpha_4^{(5)}/\alpha_5^{(5)}/}
   \plot 51 23
   \plotfmla{\fourcrit\zeta_1/\zeta_2/\alpha_5^{(5)}/\alpha_6^{(5)}/}
   \plot 51 13
   \plotfmla{\fourcrit\zeta_1/\zeta_2/\alpha_6^{(5)}/\alpha_7^{(5)}/}
   \plot 51 3
   \plotfmla{\fourcrit\zeta_1/\zeta_2/\alpha_7^{(5)}/
         \alpha_8^{(5)}{=}\alpha_3^{(6)}/}
   \plot 51 -7
   \plotfmla{\threecrit\zeta_3/\alpha_0^{(6)}{=}\alpha_0^{(7)}/
         \alpha_1^{(6)}{=}\alpha_1^{(7)}/}
   \plot 63 63
   \plotfmla{\threecrit\zeta_3/\alpha_1^{(6)}/\alpha_2^{(6)}/}
   \plot 63 53
   \plotfmla{\threecrit\zeta_3/\alpha_2^{(6)}/\alpha_3^{(6)}/}
   \plot 63 43
   \plotfmla{\threecrit\zeta_3/\alpha_3^{(6)}/\alpha_4^{(6)}/}
   \plot 63 33
   \plotfmla{\threecrit\zeta_3/\alpha_4^{(6)}/
         \alpha_5^{(6)}{=}\alpha_2^{(7)}/}
   \plot 63 23
   \plotfmla{\threecrit\zeta_3/\alpha_5^{(6)}/\alpha_6^{(6)}/}
   \plot 63 13
   \plotfmla{\threecrit\zeta_3/\alpha_6^{(6)}/\alpha_7^{(6)}/}
   \plot 63 3
   \plotfmla{\threecrit\zeta_3/\alpha_7^{(6)}/\alpha_8^{(6)}/}
   \plot 63 -7
   \vskip 1truein
\plotfigend
\vfill
\botcaption{Figure 1}
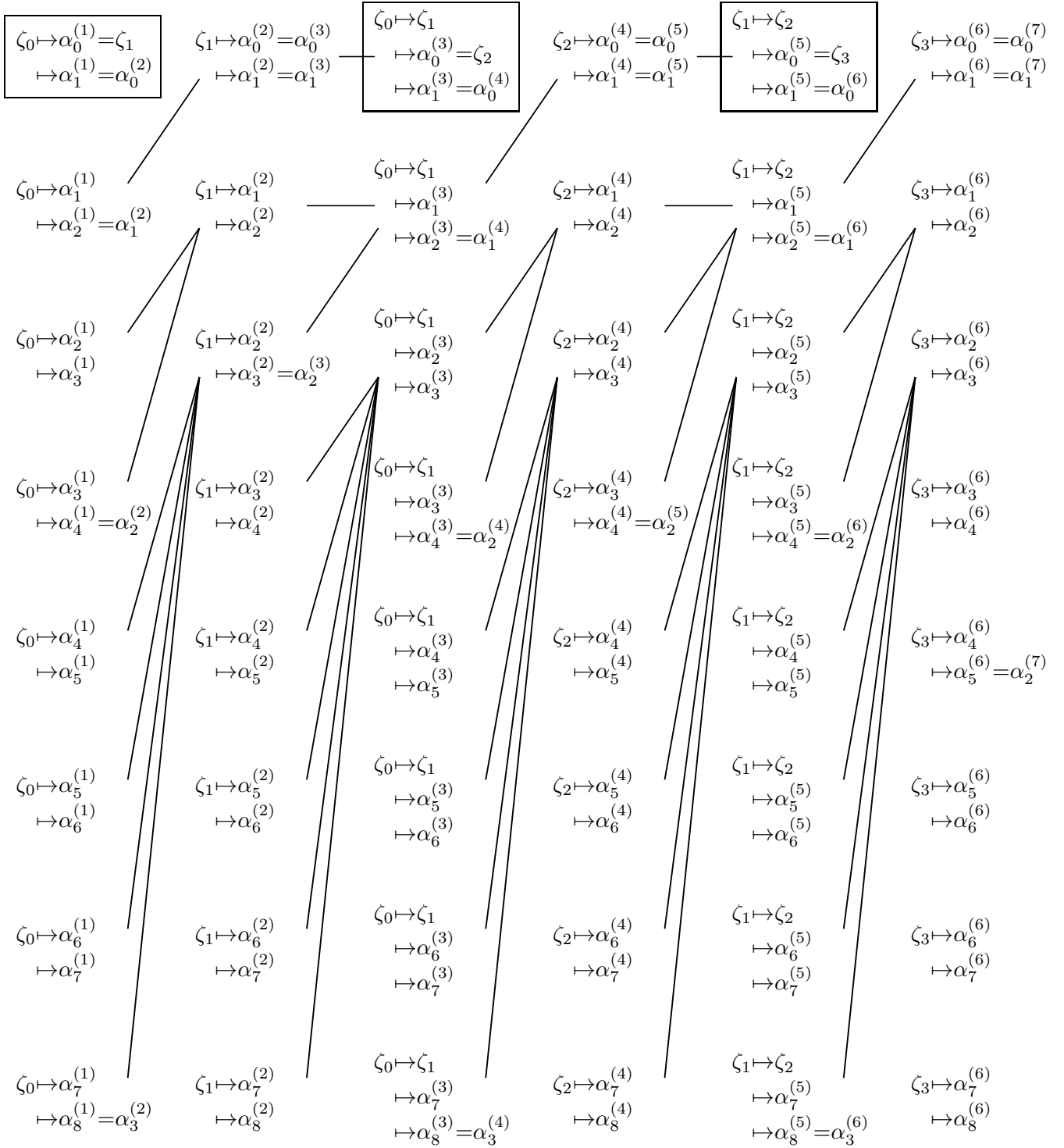
The embeddings $a_i^{(N)}$ for $1 \le N \le 6$, $0 \le i \le 7$.
\endcaption
\endinsert

The reason for making these notational changes is so that one can think of
Lemmas~3--6 as producing individual columns in a large two-dimensional
array of embeddings, as shown in Figure~1.  In this figure, each embedding
is represented by giving the relevant part of its critical sequence.
The boxed embeddings at the top of certain columns are the embeddings
$z_N$ which have to be given in advance; $z_N$ appears at the top of
column $2N+3$, and $z_0$ also appears at the top of column~$1$.
The unboxed embeddings in column $N$ are obtained by applying the
embeddings in column $N+1$ to embeddings that have already been
constructed; the lines drawn in the diagram connect each embedding in
column $N+1$ to the embeddings it is used to produce in column~$N$.

One can now state Lemmas~3--6 in a unified form, as follows:  If one is
given the first $n$ embeddings in column $N+1$, as well as all of the
boxed embeddings in preceding columns, then one can produce $C_N(n)$
embeddings in column $N$, ending up at the same critical point.
Proposition~7 just combines these column-by-column results:
given $n$ embeddings in column~$2N+3$, one can produce $C_{2N+2}(n)$
embeddings in column~$2N+\nobreak 2$, $C_{2N+1}(C_{2N+2}(n))$ embeddings in
column~$2N+1$,
and so on, leading to $C_1(C_2(\dots(C_{2N+2}(n))\dots))$
embeddings in column $1$ and hence, by Lemma~2,
$C_0(C_1(\dots(C_{2N+2}(n))\dots))$ critical points.

The process of constructing new embeddings $a_i^{(N)}$ from old ones depends
somewhat on the column number $N$.  If $N$ is odd, then the process
is relatively simple and self-contained: given the first $2^n$ embeddings
in column $N$, just apply $a_n^{(N+1)}$ to them to get the next
$2^n$ embeddings.  For $N=2$, there is a little more work: given
$C_2(n)$ embeddings in column~$2$, one must first run through the iteration
for column~$1$ to produce $2^{C_2(n)}$ embeddings there, and then apply
$a_n^{(3)}$ to the embeddings in column~$1$ to get $2^{C_2(n)}$ new embeddings
in column~$2$.  Finally, for column $2N$ where $N>1$, one has to take
two steps backward: given $C_{2N}(n)$ embeddings in column $2N$,
one must produce the corresponding $2^{C_{2N}(n)}$ embeddings in column
$2N-1$ and $C_{2N-2}(2^{C_{2N}(n)})$ embeddings in column $2N-2$ before
applying $a_n^{(2N+1)}$ to get new embeddings in column~$2N$.
In this case, though, the process of computing embeddings in column $2N-2$
from embeddings in column $2N-1$ requires going back to columns $2N-3$
and $2N-4$, and so on; in fact, one must do computations all the way
down to column $1$ in order to produce the new embeddings for column $2N$.

One can think of the entire construction of Figure~1 as being a set
of recursive equations for the embeddings~$a_i^{(N)}$:
$$a_0^{(1)} = z_0; \qquad a_0^{(2N+3)} = z_N;$$
if $C_N(n) \le i < C_N(n+1)$, then
$$a_i^{(N)} = a_n^{(N+1)}(a_{\bar \imath}^{(\hat N)}),$$
where
$$\align
\bar \imath &= i - C_N(n), \\
\hat N &= \cases N &\text{if $N$ is odd,} \\
N-1 &\text{if $N=2$,} \\
N-2 &\text{if $N>2$ is even.} \endcases
\endalign$$
This is a complicated recursion, where $a_i^{(N)}$ depends on embeddings in
both earlier and later columns; hence, if one were to take this as the
definition of the embeddings~$a_i^{(N)}$, then one would have to
check carefully that the
recursion terminates after a finite number of steps to give a definition
for $a_i^{(N)}$, and that the embeddings fit together in the proper way.
This can be done, but it is rather messy, because the order in which
the pairs $(N,i)$ are processed is complicated.  Perhaps the
simplest description of this order is in
computer science terms: the lines in Figure~1 define a collection
of trees, one rooted at each boxed entry, and the recursion proceeds
through a preorder traversal of the trees.

One final note:
We have now shown that the number $F(n)$ of critical points below
$\kappa_n$ grows very rapidly, but this does not immediately imply that
the number of critical points between $\kappa_n$ and $\kappa_{n+1}$
(which is $F(n+1)-F(n)-1$) is large for {\sl all} large $n$; it is
conceivable that $F$ grows in spurts.  (The construction in the next
section makes this seem quite plausible.)  This is not the case,
however; the preceding construction shows that there are at least
$g(n+1) - g(n) - 1$ critical points between $\zeta_n$ and
$\zeta_{n+1}$ for $n \ge 2$.  In fact, for the case $\zeta_n =
\kappa_n$, if we let $h(n) = C_1(C_2(\dots(C_{2n-4}((1))\dots))$, then
applying $a_{h(n)}^{(1)}$ to the $F(n)$ critical points below
$\kappa_n$ gives $F(n)$ additional critical points below
$\alpha_{h(n)+1}^{(1)}$; we can now apply $a_{h(n)+1}^{(1)}$ to the
current $2F(n)$ critical points to get $2F(n)$ more, and so on.
(Here we are using the fact that all of the critical points below
$\kappa_n$ are in the interval $[\kappa_0,\kappa_n)$; this is true because,
using the formula $\crit{i_1(i_2)} = i_1(\crit{i_2})$,
one can show by an easy induction on $k$ that $\crit k \ge \kappa_0$
for all $k \in \jalg$ and hence for all $k \in \jalgp$.)  This
shows that in fact $F(n+1) \ge 2^{h(n+1) - h(n)}F(n)$ for $n \ge 2$,
so $F(n+1)-F(n)$
is much larger than~$F(n)$.

\head 3. Using vectors of ordinals in the main construction \endhead

Theorem~1 and the growth rate estimate $F(n) > F_\omega(n-1)$ were
obtained by applying Corollary~8 in the case where $\zeta_n =
\kappa_n$.  If we can find a sequence of critical points $\zeta_n$
growing more slowly than~$\kappa_n$ such that embeddings as described
in Corollary~8 exist, then we get an improved lower bound on the growth
rate of $F(n)$.  In this section, we will obtain such sequences by
another application of the same construction.  We will again need the
functions $C_N$ and $g$ from the preceding section.

The goal is to produce critical points $\zeta_i$ and embeddings $z_i$
such that \fourseq z_i/\zeta_i/\zeta_{i+1}/\zeta_{i+2}/\zeta_{i+3}/.
One can compare this with Lemma~3, in which one produces embeddings
$a_i$ and critical points $\alpha_i$ such that
\threeseq a_i/\zeta/\alpha_i/\alpha_{i+1}/.  The main difference here
does not seem to be the fixed ordinal $\zeta$; after all, Lemma~5 differs
from Lemma~3 only in having two fixed ordinals instead of one, and it
is easy to prove a variant which has no fixed ordinal at all (in fact,
this will look very much like Lemma~2).  Rather, the main new difficulty
is that successive embeddings must share three varying ordinals rather
than just one.  It turns out that one can get around this difficulty by
considering ordinals in triples rather than individually.

More generally, we will consider finite sequences (`vectors') of
ordinals, which we will always assume are nonempty and in strictly
increasing order.  If $\vec\alpha$ and $\vec\beta$ are two such
sequences, $\vec\alpha < \vec\beta$ means that the last member of
$\vec\alpha$ is less than the first member of $\vec\beta$.  If
$\vec\alpha$ is such a vector of length~$L$, then its members will be
denoted $\vec\alpha(0),\vec\alpha(1),\dots, \vec\alpha(L-1)$.

Given an embedding $b$, say that {\it \fourseqv b/\vec\beta_0/
\vec\beta_1/\dots/\vec\beta_n/} if $\vec\beta_0,\dots,\vec\beta_n$ are
finite increasing sequences of ordinals of the same length,
$\vec\beta_0 < \vec\beta_1 < \dots < \vec\beta_n$, the critical point
of $b$ is $\vec\beta_0(0)$, and $b(\vec\beta_k) = \vec\beta_{k+1}$ for
each $k<n$.  (One can think of ``$b(\vec\beta_k) = \vec\beta_{k+1}$''
as an abbreviation for ``$b(\vec\beta_k(i)) = \vec\beta_{k+1}(i)$ for
all $i$'' if one wants to apply embeddings only to ordinals and other
embeddings.)  Unlike ordinary critical sequences, a critical vector
sequence is not uniquely defined; the ordinals $\vec\beta_0(i)$ for $i
\ge 1$ can be arbitrary ordinals between $\vec\beta_0(0)$ and
$\vec\beta_1(0)$.  In the applications below, these additional ordinals
will also be critical points.

\proclaim{Lemma 9} Suppose that one has vectors $\vec\tau < \vec\alpha_0
< \vec\alpha_1 < \dots < \vec\alpha_n$ of length $3$ and embeddings
$a_i \in \jalgp$ for $i < n$ such that \threeseqv
a_i/\vec\tau/\vec\alpha_i /\vec\alpha_{i+1}/ for each $i$.  Also suppose
that there exist three embeddings $z_0,z_1,z_2$ which have critical
sequences beginning $\vec\tau(0)\mapsto \vec\tau(1)\mapsto
\vec\tau(2)\mapsto \vec\alpha_0(0)$, $\vec\tau(1)\mapsto
\vec\tau(2)\mapsto \vec\alpha_0(0)\mapsto \vec\alpha_0(1)$, and
$\vec\tau(2)\mapsto \vec\alpha_0(0)\mapsto \vec\alpha_0(1)\mapsto
\vec\alpha_0(2)$, respectively.  Then there exist critical points
$\zeta_i$ ($i < 3(2^n+1)$) and embeddings $z_i$ ($i < 3 \cdot 2^n$)
such that \fourseq z_i/\zeta_i/
\zeta_{i+1}/\zeta_{i+2}/\zeta_{i+3}/ for each $i$, $\langle
\zeta_0,\zeta_1,\zeta_2\rangle = \vec\tau$, and $\langle
\zeta_{3\cdot2^n}, \zeta_{3\cdot2^n+1},
\zeta_{3\cdot2^n+2}\rangle = \vec\alpha_n$.\endproclaim

\demo{Proof} Induct on $n$; for $n=0$, just let $z_i$ be as given for
$i=0,1,2$, and let $\zeta_i = \vec\tau(i)$, $\zeta_{3+i} =
\vec\alpha_0(i)$.  Given the result for $n$, if one has $\vec\alpha_i$ for $i
\le n+1$ and $a_i$ for $i < n+1$ as above, then the induction
hypothesis gives $\zeta_i$ ($i < 3(2^n+1)$) and $z_i$ ($i < 3\cdot
2^n$).  Let $\zeta_{3\cdot 2^n+i} = a_n(\zeta_i)$ for $3 \le i <
3\cdot2^n+3$; this equation holds also for $i=0,1,2$, since
$a_n(\vec\tau) = \vec\alpha_n$.  Hence, if we let $z_{3\cdot 2^n+i} =
a_n(z_i)$ for $i < 3\cdot2^n$, then the embeddings $z_{3\cdot 2^n+i}$
will have the required critical sequences.  Finally, $\langle
\zeta_{3\cdot2^{n+1}}, \zeta_{3\cdot2^{n+1}+1},
\zeta_{3\cdot2^{n+1}+2}\rangle = a_n(\vec\alpha_n) = \vec\alpha_{n+1}$.
This completes the induction. \QED

We now have reduced the problem of obtaining embeddings $z_i$ related
as above to that of obtaining embeddings $a_i$ which map triples of
ordinals as above.  But it is not hard to see that the methods used
for Lemmas~3--6 work without change if one uses vectors of ordinals
(of a fixed finite length) instead of individual ordinals.  For instance,
the `vectorized' version of Lemma~3 is:

\proclaim{Lemma 10} Suppose that one has vectors $\vec\tau < \vec\tau' <
\vec\beta_0 < \vec\beta_1 < \dots < \vec\beta_n$ and embeddings $b_i
\in \jalgp$ for $i < n$ such that \threeseqv b_i/\vec\tau'/\vec\beta_i
/\vec\beta_{i+1}/ for each $i$.  Also suppose that there is an
embedding $t$ such that \threeseqv t/\vec\tau/\vec\tau' /\vec\beta_0/.
Then there exist vectors $\vec\alpha_i$ ($i \le 2^n$) and embeddings
$a_i$ ($i < 2^n$) such that \threeseqv a_i/\vec\tau/\vec\alpha_i
/\vec\alpha_{i+1}/ for each $i$, $\vec\alpha_0 = \vec\tau'$, and
$\vec\alpha_{2^n} = \vec\beta_n$.\endproclaim

The proof of this lemma is the same as that of the original Lemma~3.
(Here we have $b_n(\vec\tau) = \vec\tau$ because all of the ordinals
$\vec\tau(i)$ are below $\vec\tau'(0)$, which is the critical point
of~$b_n$.)  It is just as easy to produce vector versions of Lemmas~4--6.
One can now combine these as before to get the following analogue of
Corollary~8:

\proclaim{Lemma 11} Suppose one has vectors
$\vec\tau_0 < \vec\tau_1 < \dots < \vec\tau_m$ of length\/ $3$ and embeddings
$t_k$ ($0 \le k \le m-3$) such that \fourseqv t_k/ \vec\tau_k/
\vec\tau_{k+1}/ \vec\tau_{k+2}/ \vec\tau_{k+3}/.  Suppose also that
there exist three embeddings $z_0,z_1,z_2$ which have critical
sequences beginning $\vec\tau_0(0)\mapsto \vec\tau_0(1)\mapsto
\vec\tau_0(2)\mapsto \vec\tau_1(0)$, $\vec\tau_0(1)\mapsto
\vec\tau_0(2)\mapsto \vec\tau_1(0)\mapsto \vec\tau_1(1)$, and
$\vec\tau_0(2)\mapsto \vec\tau_1(0)\mapsto \vec\tau_1(1)\mapsto
\vec\tau_1(2)$, respectively.  
Then there exist critical points
$\zeta_i$ ($i < 3(g(n)+1)$) and embeddings $z_i$ ($i < 3 g(n)$)
such that \fourseq z_i/\zeta_i/
\zeta_{i+1}/\zeta_{i+2}/\zeta_{i+3}/ for each $i$, $\langle
\zeta_0,\zeta_1,\zeta_2\rangle = \vec\tau_0$, and $\langle
\zeta_{3g(m)}, \zeta_{3g(m)+1},
\zeta_{3g(m)+2}\rangle = \vec\tau_m$.\endproclaim

It follows that, if the hypotheses of Lemma~11 hold, then the number of
critical points below $\vec\tau_m(0)$ is at least
$g(3g(m))$.  For instance, if we let $\vec\tau_k = \langle \kappa_{3k},
\kappa_{3k+1}, \kappa_{3k+2} \rangle$, $t_k = j^{[3k]} \circ j^{[3k]}
\circ j^{[3k]}$, and $z_i = j^{[i]}$ for $i=0,1,2$, then we get
$F(3m) \ge g(3g(m))$, which is a much better bound than~$g(3m)$
for $m \ge 3$.

The hypotheses of Lemma~11 look very similar to those of Corollary~8;
one merely has to have suitable embeddings working on triples of
ordinals rather than ordinals.  Hence, in order to meet these
hypotheses, we merely have to repeat the modifications above,
replacing triples of ordinals with triples of triples of ordinals,
or, equivalently, vectors of length~$9$.  More generally, the proof
of Lemma~11 goes through without change to produce the following
version:

\proclaim{Lemma 12} Suppose one has vectors $\vec\tau_0 < \vec\tau_1 <
\dots < \vec\tau_m$ of length\/ $3L$ and embeddings $t_k$ ($0 \le k \le
m-3$) such that \fourseqv t_k/ \vec\tau_k/ \vec\tau_{k+1}/ \vec\tau_{k+2}/
\vec\tau_{k+3}/.  Let $\vec\tau_k^{[i]}$ for $i=0,1,2$ be the three
sequences of length $L$ composing $\vec\tau_k$ (i.e., $\vec\tau_k =
\vec\tau_k^{[0]}{}^\cap \vec\tau_k^{[1]}{}^\cap \vec\tau_k^{[2]}$), and
suppose that there exist three embeddings $z_0,z_1,z_2$ which have
critical vector sequences beginning $\vec\tau_0^{[0]}\mapsto
\vec\tau_0^{[1]}\mapsto \vec\tau_0^{[2]}\mapsto \vec\tau_1^{[0]}$,
$\vec\tau_0^{[1]}\mapsto \vec\tau_0^{[2]}\mapsto \vec\tau_1^{[0]}\mapsto
\vec\tau_1^{[1]}$, and $\vec\tau_0^{[2]}\mapsto \vec\tau_1^{[0]}\mapsto
\vec\tau_1^{[1]}\mapsto \vec\tau_1^{[2]}$, respectively.  Then there exist
vectors $\vec\zeta_i$ ($i < 3(g(n)+1)$) and embeddings $z_i$ ($i < 3
g(n)$) such that \fourseqv z_i/\vec\zeta_i/
\vec\zeta_{i+1}/\vec\zeta_{i+2}/\vec\zeta_{i+3}/ for each $i$,
$\vec\zeta_0{}^\cap \vec\zeta_1{}^\cap \vec\zeta_2 =
\vec\tau_0$, and $\vec\zeta_{3g(m)}{}^\cap
\vec\zeta_{3g(m)+1}{}^\cap \vec\zeta_{3g(m)+2} =
\vec\tau_m$.\endproclaim

Actually, we have $\vec\zeta_{3g(k)}{}^\cap
\vec\zeta_{3g(k)+1}{}^\cap \vec\zeta_{3g(k)+2} =
\vec\tau_k$ for all $k \le m$.

By using Lemma~12 repeatedly, we can generate critical points according
to iterates of the function $g$.  Let $g_m(n) = g((3g)^m(n))$; that is,
$g_0(n) = g(n)$ and $g_{m+1}(n) = g_m(3g(n))$.  (One also has
$g_{m+1}(n) = g(3g_m(n))$.)  For any fixed $m$, a suitable
starting point for applying Lemma~12 is the sequences $\vec\tau_n =
\langle \kappa_{3^mn+i} \colon 0 \le i < 3^m \rangle$, since we can define
$t_n$ to be $j^{3^mn}(j^{3^m})$, where $j^l$ is the composition of $l$ $j$'s.
The three additional embeddings $z_i$ are just $j^{3^{m-1}i}(j^{3^{m-1}})$,
$i=0,1,2$.  Lemma~12 now lets us define $z_n$ and $\vec\zeta_n$ (of
length $3^{m-1}$) for all $n$, and we have $\vec\zeta_{3g(n)}(0) =
\vec\tau_n(0) = \kappa_{3^mn}$.  Furthermore, we have $\vec\zeta_n =
\langle \kappa_{3^{m-1}n+i} \colon 0 \le i < 3^{m-1} \rangle$ for $n \le 5$
(for $n \le 8$, actually, but this is irrelevant), so we can define three
embeddings $z'_i = j^{3^{m-2}i}(j^{3^{m-2}})$, $i = 0,1,2$, and again we are
in a position to apply Lemma~12.  Continuing this through a total of $m-1$
applications of Lemma~12 leaves us with ordinal sequences of length $3$,
and we are now in a position to apply Lemma~11 and then
Corollary~8.  We therefore get the following result: for all $m$ and $n$,
$F(3^mn) \ge g_m(n)$.

In particular, we have $F(3^m) \ge g_m(1)$ for all $m$.  Now, a direct
computation shows that $g_1(1) = 4 > 3 = F_{\omega+1}(0)$.  This, together
with the equation $g_{m+1}(n) = g(3g_m(n))$ and the inequality
$g(n) > F_\omega(n-1)$ for $n \ge 4$ from the preceding section,
yields an easy inductive proof that $g_m(1) > F_{\omega+1}(m-1)$ for all
$m \ge 1$.  So $F(3^m) > F_{\omega+1}(m-1)$.  If we apply this in the case
where $m = \lfloor \log_3(n) \rfloor$, we get the result stated in the
introduction:  $F(n) > F_{\omega+1}(\lfloor \log_3(n) \rfloor-1)$ for
all $n \ge 3$.

%

%
\head 4. Irregular constructions \endhead

The preceding sections have shown that $F(n)$, the number of critical points
below $\kappa_n$, is very large when $n$ is large.  Now we will turn to the
specific value $F(4)$.  (This is the first value of $F$ that could be large;
we will note in the next section that $F(3) = 4$.)  The lower bound for
$F(4)$ given by section~2 is $g(4) = 256$, and section~3 does
not improve this value.  Here we will show, by use of some {\it ad hoc}
modifications of the main construction, that $F(4)$ is quite large.

The numerical functions involved will be somewhat different from those
used in section~2; the new functions will be called $\Ct_N$ rather than
$C_N$.  Also, the resulting critical points and embeddings will be referred
to as $\at_i^{(N)}$ and $\et_i^{(N)}$ rather than $\alpha_i^{(N)}$ and $a_i^{(N)}$ as at the end of section~2.

The first four stages of the construction, Lemmas~2--5 (or the first
four columns of Figure~1), will remain unchanged: $\Ct_N$ is the same
as $C_N$ for $N \le 3$.  Since this is a specific construction, we will
go ahead and use the specific values
$\kappa_0,\kappa_1,\kappa_2,\kappa_3$ rather than the general
$\zeta_0,\zeta_1,\zeta_2,\zeta_3$.  In other words, we will use
Lemma~5 with $\zeta=\kappa_0$, $\zeta'=\kappa_1$, $\zeta''=
\kappa_2$, and $z=j$.

As a first example, suppose that we have an embedding with critical
point $\kappa_2$ which sends $\kappa_2$ to $\kappa_3$ and $\kappa_3$ to
some critical point $\mu$.  Then we can apply Lemma~5 with $n=1$,
followed by Lemmas~4,~3, and~2, to produce $2^{2^{C_2(2^1)}} = 256$
critical points below $\mu$.  Hence, if $\mu$ is less than $\kappa_4$,
then $F(4) > 256$.

We will see in the next section that such an embedding exists, with $\mu <
\kappa_4$; in fact, $\jsub{12}$ will work, where $\jsub n$ is as defined in
section~1.  For the rest of this section, we will continue to construct many
critical points below $\kappa_4$, assuming that we can get certain specific
embeddings; the proofs that such embeddings exist will be given in the next
section.

For the next step, we will construct two additional stages (columns) by
a method different from that used in section~2.  First, we give the recursive
definitions of the relevant functions $\Ct_4$ and $\Ct_5$:
$$\gather \Ct_4(0)=0,\qquad \Ct_4(1)=8,\qquad \Ct_4(m+1)=\Ct_4(m) +
\Ct_1(\Ct_2(\Ct_3(\Ct_4(m))))-1\quad\text{for $m>0$;} \\
\Ct_5(0)=0,\qquad \Ct_5(1)=2,\qquad \Ct_5(m+1)=\Ct_5(m) +
\Ct_3(\Ct_4(\Ct_5(m)))\quad\text{for $m>0$.} \endgather$$
At this point, it will be useful to introduce an abbreviation: for
$N < M$, let
$$\Ctfunc NM(n) = \Ct_N(\Ct_{N+1}(\dots(\Ct_{M-1}(n))\dots)).$$
So the recursive formula for $\Ct_4$ can be written as
$\Ct_4(m+1) = \Ct_4(m) + \Ctfunc14(\Ct_4(m))-1$.  (One could
further abbreviate this to $\Ct_4(m) + \Ctfunc15(m)-1$, but it
seems preferable to keep the function being defined explicit.)
And now the corresponding lemmas:

\proclaim{Lemma 13} Suppose that one has critical points $\at^{(5)}_0 <
\at^{(5)}_1 < \dots < \at^{(5)}_n$, where $n > 0$, and embeddings
$\et^{(5)}_i$ for $i < n$ such that \fourseq
\et^{(5)}_i/\kappa_0/\kappa_2/\at^{(5)}_i /\at^{(5)}_{i+1}/ for each
$i$.  Also suppose that $\et^{(5)}_0(\kappa_1) = \kappa_3$.  Then
there exist critical points $\at^{(4)}_i$ ($i \le \Ct_4(n)$) and
embeddings $\et^{(4)}_i$ ($i < \Ct_4(n)$) such that \threeseq
\et^{(4)}_i/ \kappa_2/ \at^{(4)}_i / \at^{(4)}_{i+1}/ for each $i$,
$\at^{(4)}_0 = \kappa_3$, $\at^{(4)}_1 = \at^{(5)}_0$, and
$\at^{(4)}_{\Ct_4(n)} = \at^{(5)}_n$.\endproclaim

\demo{Proof} Induct on $n \ge 1$.  For $n=1$, first let $\et^{(4)}_0 =
\et^{(5)}_0(j)$, $\at^{(4)}_0 = \kappa_3$, and $\at^{(4)}_1 =
\at^{(5)}_0$.  Now apply Lemma~5 (with $z=j$), Lemma~4, and Lemma~3
(with $z=j$ again) to get $\Ct_1(\Ct_2(\Ct_3(1))) = 8$ embeddings
$\et^{(1)}_i$ and corresponding critical points $\at^{(1)}_i$ ($i \le
8$) such that \threeseq
\et^{(1)}_i/\kappa_0/\at^{(1)}_i/\at^{(1)}_{i+1}/, and $\at^{(1)}_8 =
\at^{(5)}_0$; since we used $z=j$ in Lemma~3, we will have $\et^{(1)}_0
= j$, $\at^{(1)}_0 = \kappa_1$, and $\at^{(1)}_1 = \kappa_2$.  Now let
$\et^{(4)}_i = \et^{(5)}_0(\et^{(1)}_i)$ and $\at^{(4)}_{i+1} =
\et^{(5)}_0(\at^{(1)}_{i+1})$ for $1 \le i < 8$; these equations also
hold for $i=0$, and we get $\at^{(4)}_8 = \at^{(5)}_1$, so everything
fits together and the case $n=1$ is done.

The induction step is similar.  Given the result for $n$, if one has
$\at^{(5)}_i$ for $i \le n+1$ and $\et^{(5)}_i$ for $i < n+1$ as above,
then the induction hypothesis gives critical points $\at^{(4)}_i$ ($i
\le \Ct_4(n)$) and embeddings $\et^{(4)}_i$ ($i < \Ct_4(n)$) as
required, with $\at^{(4)}_0 = \kappa_3$, $\at^{(4)}_1 = \at^{(5)}_0$,
and $\at^{(4)}_{\Ct_4(n)} = \at^{(5)}_n$.  We can now apply
Lemmas~5,~4, and~3 successively to get embeddings $\et^{(1)}_i$ ($i <
\Ctfunc14(\Ct_4(n))$) and critical points $\at^{(1)}_i$ ($i \le
\Ctfunc14(\Ct_4(n))$) such that \threeseq
\et^{(1)}_i/\kappa_0/\at^{(1)}_i/\at^{(1)}_{i+1}/; also, we will have
$\at^{(1)}_1 = \kappa_2$ and $\at^{(1)}_{\Ctfunc14(\Ct_4(n))} =
\at^{(5)}_n$.  Now let $\et^{(4)}_{\Ct_4(n)-1+i} =
\et^{(5)}_n(\et^{(1)}_i)$ and $\at^{(4)}_{\Ct_4(n)+i} =
\et^{(5)}_n(\at^{(1)}_{i+1})$ for $1 \le i < \Ctfunc14(\Ct_4(n))$; this
fits together suitably because $\et^{(5)}_n(\at^{(1)}_1) = \at^{(5)}_n
= \at^{(4)}_{\Ct_4(n)}$, and we also get
$\at^{(4)}_{\Ct_4(n+1)} = \at^{(5)}_{n+1}$, so the induction is complete. \QED

\proclaim{Lemma 14} Suppose that one has critical points $\at^{(6)}_0 <
\at^{(6)}_1 < \dots < \at^{(6)}_n$, where $n > 0$, and embeddings
$\et^{(6)}_i$ for $i < n$ such that \fourseq
\et^{(6)}_i/\kappa_1/\kappa_2/\at^{(6)}_i /\at^{(6)}_{i+1}/ for each
$i$.  Also suppose that $\et^{(6)}_0(j)(\kappa_1) = \kappa_3$.  Then
there exist critical points $\at^{(5)}_i$ ($i \le \Ct_5(n)$) and
embeddings $\et^{(5)}_i$ ($i < \Ct_5(n)$) such that \fourseq
\et^{(5)}_i/\kappa_0/\kappa_2/\at^{(5)}_i /\at^{(5)}_{i+1}/ for each $i$,
$\et^{(5)}_0(\kappa_1) = \kappa_3$, $\at^{(5)}_0 = \at^{(6)}_0$, and
$\at^{(5)}_{\Ct_5(n)} = \at^{(6)}_n$.\endproclaim

\demo{Proof} Induct on $n$.  For $n=1$, let $\et^{(5)}_0 =
\et^{(6)}_0(j)$, $\et^{(5)}_1 = \et^{(6)}_0(\et^{(5)}_0(j)(j))$,
$\at^{(5)}_0 = \at^{(6)}_0$, $\at^{(5)}_1 = \et^{(6)}_0(\kappa_3)$, and
$\at^{(5)}_2 = \at^{(6)}_1$.  [This is just a shorter way of stating
the following: Let $\et^{(5)}_0 = \et^{(6)}_0(j)$, $\at^{(5)}_0 =
\at^{(6)}_0$, and $\at^{(5)}_1 = \et^{(6)}_0(\kappa_3)$.  Apply
Lemma~13 with $n=1$, and then apply Lemma~5 with $n=1$ (not $n=8$), to
get embeddings $\et^{(3)}_i$ ($i < 2$) and critical points
$\at^{(3)}_i$ ($i \le 2$) such that \fourseq \et^{(3)}_i/ \kappa_0/
\kappa_1/ \at^{(3)}_i/ \at^{(3)}_{i+1}/, and $\at^{(3)}_0 = \kappa_2$,
$\at^{(3)}_1 = \kappa_3$, and $\at^{(3)}_2 = \at^{(6)}_0$.  Now let
$\et^{(5)}_1 = \et^{(6)}_0(\et^{(3)}_1)$ and $\at^{(5)}_2 =
\at^{(6)}_1$.]

Given the result for $n$, if one has $\at^{(6)}_i$ for $i \le n+1$ and
$\et^{(6)}_i$ for $i < n+1$ as above, then the induction hypothesis
gives critical points $\at^{(5)}_i$ ($i \le \Ct_5(n)$) and embeddings
$\et^{(5)}_i$ ($i < \Ct_5(n)$) as required, with $\et^{(5)}_0(\kappa_1)
= \kappa_3$, $\at^{(5)}_0 = \at^{(6)}_0$, and $\at^{(5)}_{\Ct_5(n)} =
\at^{(6)}_n$.  We can now apply Lemmas~13 and~5 successively to get
embeddings $\et^{(3)}_i$ ($i < \Ctfunc35(\Ct_5(n))$) and critical
points $\at^{(3)}_i$ ($i \le \Ctfunc35(\Ct_5(n))$) such that \fourseq
\et^{(3)}_i/\kappa_0/\kappa_1/\at^{(3)}_i/\at^{(3)}_{i+1}/; also, also,
we will have $\at^{(3)}_0 = \kappa_2$ and
$\at^{(3)}_{\Ctfunc35(\Ct_5(n))} = \at^{(6)}_n$.  Now let
$\et^{(5)}_{\Ct_5(n)+i} = \et^{(6)}_n(\et^{(3)}_i)$ and
$\at^{(5)}_{\Ct_5(n)+i+1} = \et^{(6)}_n(\at^{(3)}_{i+1})$ for $0 \le i
< \Ctfunc35(\Ct_5(n))$; this fits together suitably because
$\et^{(6)}_n(\at^{(3)}_0) = \at^{(6)}_n = \at^{(5)}_{\Ct_5(n)}$, and we
also get $\at^{(5)}_{\Ct_5(n+1)} = \at^{(6)}_{n+1}$, so the induction
is complete. \QED

The modified construction to this point is shown in Figure~2.  (For
now, ignore the references to $\at^{(7)}_i$.)  The main difference
between this and the original construction is that here we extend
column~4 by applying embeddings in column 5 to those in column 1
instead of column~2, and we extend column~5 by applying embeddings in
column 6 to those in column 3 instead of column~5.

\topinsert
\plotfigbegin
   \normalbaselines
   \vskip 5.7truein
   \plotPSbegin
      \setdoline
      gsave
      -9 53 translate
      1 1 0 doline
      1 2 1 doline
      1 3 1 doline
      1 4 2 doline
      1 5 2 doline
      2 -1 0 doline
      2 1 1 doline
      2 2 1 doline
      2 3 2 doline
      2 4 2 doline
      2 5 2 doline
      3 1 0 doline
      3 2 1 doline
      3 3 1 doline
      3 4 2 doline
      3 5 2 doline
      4 -1 0 doline
      4 1 0 doline
      4 2 0 doline
      4 3 0 doline
      4 4 0 doline
      4 5 0 doline
      5 -1 0 doline
      5 1 0 doline
      5 2 1 doline
      5 3 1 doline
      5 4 1 doline
      5 5 1 doline
      grestore
      0.4 72 div setlinewidth
      [1 0.5] 0 setdash
      newpath
      69 48 moveto
      57.98 48 lineto
      56.02 38 lineto
      45.98 38 lineto
      45 13 lineto
      45 0 lineto
      stroke
   \plotPSend
   \plotfmla{\boxed{\threecrit\kappa_0/\at_0^{(1)}{=}\kappa_1/
         \at_1^{(1)}{=}\at_0^{(2)}/}}
   \plot 3 53
   \plotfmla{\threecrit\kappa_0/\at_1^{(1)}/\at_2^{(1)}{=}\at_1^{(2)}/}
   \plot 3 43
   \plotfmla{\threecrit\kappa_0/\at_2^{(1)}/\at_3^{(1)}/}
   \plot 3 33
   \plotfmla{\threecrit\kappa_0/\at_3^{(1)}/\at_4^{(1)}{=}\at_2^{(2)}/}
   \plot 3 23
   \plotfmla{\threecrit\kappa_0/\at_4^{(1)}/\at_5^{(1)}/}
   \plot 3 13
   \plotfmla{\threecrit\kappa_0/\at_5^{(1)}/\at_6^{(1)}/}
   \plot 3 3
   \plotfmla{\threecrit\kappa_1/\at_0^{(2)}{=}\at_0^{(3)}/
         \at_1^{(2)}{=}\at_1^{(3)}/}
   \plot 15 53
   \plotfmla{\threecrit\kappa_1/\at_1^{(2)}/\at_2^{(2)}/}
   \plot 15 43
   \plotfmla{\threecrit\kappa_1/\at_2^{(2)}/\at_3^{(2)}{=}\at_2^{(3)}/}
   \plot 15 33
   \plotfmla{\threecrit\kappa_1/\at_3^{(2)}/\at_4^{(2)}/}
   \plot 15 23
   \plotfmla{\threecrit\kappa_1/\at_4^{(2)}/\at_5^{(2)}/}
   \plot 15 13
   \plotfmla{\threecrit\kappa_1/\at_5^{(2)}/\at_6^{(2)}/}
   \plot 15 3
   \plotfmla{\boxed{\fourcrit\kappa_0/\kappa_1/\at_0^{(3)}{=}\kappa_2/
         \at_1^{(3)}{=}\at_0^{(4)}/}}
   \plot 27 53
   \plotfmla{\fourcrit\kappa_0/\kappa_1/\at_1^{(3)}/\at_2^{(3)}{=}\at_1^{(4)}/}
   \plot 27 43
   \plotfmla{\fourcrit\kappa_0/\kappa_1/\at_2^{(3)}/\at_3^{(3)}/}
   \plot 27 33
   \plotfmla{\fourcrit\kappa_0/\kappa_1/\at_3^{(3)}/\at_4^{(3)}{=}\at_2^{(4)}/}
   \plot 27 23
   \plotfmla{\fourcrit\kappa_0/\kappa_1/\at_4^{(3)}/\at_5^{(3)}/}
   \plot 27 13
   \plotfmla{\fourcrit\kappa_0/\kappa_1/\at_5^{(3)}/\at_6^{(3)}/}
   \plot 27 3
   \plotfmla{\threecrit\kappa_2/\at_0^{(4)}{=}\kappa_3/
         \at_1^{(4)}{=}\at_0^{(5)}/}
   \plot 39 53
   \plotfmla{\threecrit\kappa_2/\at_1^{(4)}/\at_2^{(4)}/}
   \plot 39 43
   \plotfmla{\threecrit\kappa_2/\at_2^{(4)}/\at_3^{(4)}/}
   \plot 39 33
   \plotfmla{\threecrit\kappa_2/\at_3^{(4)}/\at_4^{(4)}/}
   \plot 39 23
   \plotfmla{\threecrit\kappa_2/\at_4^{(4)}/\at_5^{(4)}/}
   \plot 39 13
   \plotfmla{\threecrit\kappa_2/\at_5^{(4)}/\at_6^{(4)}/}
   \plot 39 3
   \plotfmla{\fourcrit\kappa_0/\kappa_2/\at_0^{(5)}{=}\at_0^{(6)}/\at_1^{(5)}/}
   \plot 51 53
   \plotfmla{\fourcrit\kappa_0/\kappa_2/\at_1^{(5)}/\at_2^{(5)}{=}\at_1^{(6)}/}
   \plot 51 43
   \plotfmla{\fourcrit\kappa_0/\kappa_2/\at_2^{(5)}/\at_3^{(5)}/}
   \plot 51 33
   \plotfmla{\fourcrit\kappa_0/\kappa_2/\at_3^{(5)}/\at_4^{(5)}/}
   \plot 51 23
   \plotfmla{\fourcrit\kappa_0/\kappa_2/\at_4^{(5)}/\at_5^{(5)}/}
   \plot 51 13
   \plotfmla{\fourcrit\kappa_0/\kappa_2/\at_5^{(5)}/\at_6^{(5)}/}
   \plot 51 3
   \plotfmla{\boxed{\fourcrit\kappa_1/\kappa_2/\at_0^{(6)}{=}\mu/
         \at_1^{(6)}{=}\at_0^{(7)}/}}
   \plot 63 53
   \plotfmla{\fourcrit\kappa_1/\kappa_2/\at_1^{(6)}/\at_2^{(6)}{=}\at_1^{(7)}/}
   \plot 63 43
   \plotfmla{\fourcrit\kappa_1/\kappa_2/\at_2^{(6)}/\at_3^{(6)}/}
   \plot 63 33
   \plotfmla{\fourcrit\kappa_1/\kappa_2/\at_3^{(6)}/\at_4^{(6)}{=}\at_2^{(7)}/}
   \plot 63 23
   \plotfmla{\fourcrit\kappa_1/\kappa_2/\at_4^{(6)}/\at_5^{(6)}/}
   \plot 63 13
   \plotfmla{\fourcrit\kappa_1/\kappa_2/\at_5^{(6)}/\at_6^{(6)}/}
   \plot 63 3
\plotfigend
\botcaption{Figure 2}
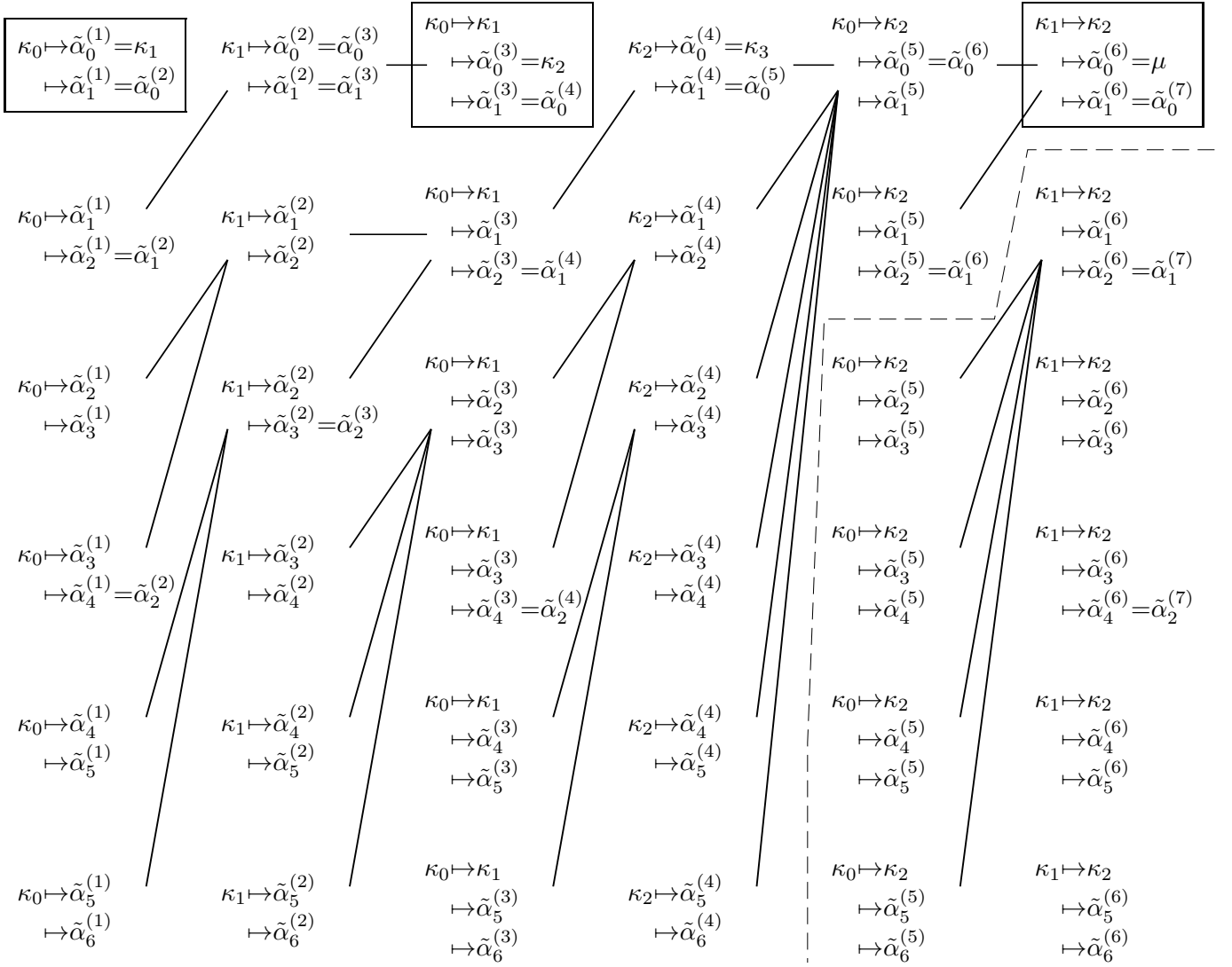
The embeddings $\et_m^{(n)}$ for $1 \le n \le 6$, $0 \le m \le 5$.
\endcaption
\endinsert

We now assume the existence of an embedding $k$ such that
\fourseq k/\kappa_1 / \kappa_2 / \mu / \nu/ for some $\mu$ and
$\nu$, satisfying the additional condition $k(j)(\kappa_1) =
\kappa_3$.  (We will see in the next section that $k = \jsub {10}$
works.)  Then we can let $\et^{(6)}_0 = k$, $\at^{(6)}_0 = \mu$, and
$\at^{(6)}_1 = \nu$, and apply Lemmas~14,~13, 5, 4, 3, and~2
successively to fill in all of Figure~2 except the part below the
dashed line; this generates a total of $\Ctfunc 06(1)$ critical points
below $\nu$.

\goodbreak

Since we have $\Ct_ 2 = C_2$, we can use the estimate $\Ct_ 2(m+4) \ge
F_4(m)+3$.  This leads to the computations $\Ct_4(1) = 8$,
$\Ctfunc35(1) = 256$, $\Ctfunc25(1) \ge F_4(252)+3$, $\Ctfunc15(1) >
F_4(253)$, $\Ct_4(2) > F_4(253)+7$, $\Ctfunc35(2) > F_4(254)+4$, and
$\Ctfunc06(1) = \Ctfunc05(2) > \Ctfunc 25(2) > F_4(F_4(254))$, so there
are more than $F_4(F_4(254))$ critical points below $\nu$.  The value
of $\nu$ produced in the next section will be less than $\kappa_4$, so
we will have $F(4) > F_4(F_4(254)) > F_5(1)$.  This already shows that
$F(4)$ is quite large.

To get even larger lower bounds on $F(4)$, we will add more
columns to the irregular construction.  To do this, we assume the existence
of two additional embeddings $k',k''$ such that \fourseq
k'/\kappa_0 / \mu / \nu / \xi/ (where $\mu$ and $\nu$
are the same as for the embedding $k$), $k''$ has critical sequence
beginning $\kappa_2 \mapsto \nu$, and $k''(k)(j)(\mu) = \xi$.  (In the next
section, we will verify that the embeddings $k' = \jsub {10}(\jsub {11})$ and
$k'' = \jsub 9(\jsub {14})$ satisfy these conditions.)

\topinsert
\plotfigbegin
   \normalbaselines
   \vskip 5.7truein
   \plotPSbegin
      \setdoline
      -8 53 translate
      1 -1 0 doline
      1 1 0 doline
      1 2 0 doline
      1 3 0 doline
      1 4 0 doline
      1 5 0 doline
      2 1 0 doline
      2 2 1 doline
      2 3 1 doline
      2 4 2 doline
      2 5 2 doline
      3 0 0 doline
      3 1 0 doline
      3 2 0 doline
      3 3 0 doline
      3 4 0 doline
      3 5 0 doline
      4 0 0 doline
      4 1 0 doline
      4 2 0 doline
      4 3 0 doline
      4 4 0 doline
      4 5 0 doline
      5 0 0 doline
      5 1 0 doline
      5 2 0 doline
      5 3 0 doline
      5 4 0 doline
      5 5 0 doline
   \plotPSend
   \plotfmla{\threecrit\mu/\at_0^{(7)}{=}\at_0^{(8)}/\at_1^{(7)}/}
   \plot 4 53
   \plotfmla{\threecrit\mu/\at_1^{(7)}/\at_2^{(7)}/}
   \plot 4 43
   \plotfmla{\threecrit\mu/\at_2^{(7)}/\at_3^{(7)}/}
   \plot 4 33
   \plotfmla{\threecrit\mu/\at_3^{(7)}/\at_4^{(7)}/}
   \plot 4 23
   \plotfmla{\threecrit\mu/\at_4^{(7)}/\at_5^{(7)}/}
   \plot 4 13
   \plotfmla{\threecrit\mu/\at_5^{(7)}/\at_6^{(7)}/}
   \plot 4 3
   \plotfmla{\boxed{\fourcrit\kappa_0/\mu/\at_0^{(8)}{=}\nu/
         \at_1^{(8)}{=}\at_0^{(9)}/}}
   \plot 16 53
   \plotfmla{\fourcrit\kappa_0/\mu/\at_1^{(8)}/\at_2^{(8)}{=}\at_1^{(9)}/}
   \plot 16 43
   \plotfmla{\fourcrit\kappa_0/\mu/\at_2^{(8)}/\at_3^{(8)}/}
   \plot 16 33
   \plotfmla{\fourcrit\kappa_0/\mu/\at_3^{(8)}/\at_4^{(8)}{=}\at_2^{(9)}/}
   \plot 16 23
   \plotfmla{\fourcrit\kappa_0/\mu/\at_4^{(8)}/\at_5^{(8)}/}
   \plot 16 13
   \plotfmla{\fourcrit\kappa_0/\mu/\at_5^{(8)}/\at_6^{(8)}/}
   \plot 16 3
   \plotfmla{\threecrit\nu/\at_0^{(9)}{=}\xi/\at_1^{(9)}/}
   \plot 28 53
   \plotfmla{\threecrit\nu/\at_1^{(9)}/\at_2^{(9)}/}
   \plot 28 43
   \plotfmla{\threecrit\nu/\at_2^{(9)}/\at_3^{(9)}/}
   \plot 28 33
   \plotfmla{\threecrit\nu/\at_3^{(9)}/\at_4^{(9)}/}
   \plot 28 23
   \plotfmla{\threecrit\nu/\at_4^{(9)}/\at_5^{(9)}/}
   \plot 28 13
   \plotfmla{\threecrit\nu/\at_5^{(9)}/\at_6^{(9)}/}
   \plot 28 3
   \plotfmla{\fourcrit\kappa_0/\nu/\at_0^{(10)}/\at_1^{(10)}/}
   \plot 40 53
   \plotfmla{\fourcrit\kappa_0/\nu/\at_1^{(10)}/\at_2^{(10)}/}
   \plot 40 43
   \plotfmla{\fourcrit\kappa_0/\nu/\at_2^{(10)}/\at_3^{(10)}/}
   \plot 40 33
   \plotfmla{\fourcrit\kappa_0/\nu/\at_3^{(10)}/\at_4^{(10)}/}
   \plot 40 23
   \plotfmla{\fourcrit\kappa_0/\nu/\at_4^{(10)}/\at_5^{(10)}/}
   \plot 40 13
   \plotfmla{\fourcrit\kappa_0/\nu/\at_5^{(10)}/\at_6^{(10)}/}
   \plot 40 3
   \plotfmla{\fourcrit\kappa_1/\nu/\at_0^{(11)}/\at_1^{(11)}/}
   \plot 52 53
   \plotfmla{\fourcrit\kappa_1/\nu/\at_1^{(11)}/\at_2^{(11)}/}
   \plot 52 43
   \plotfmla{\fourcrit\kappa_1/\nu/\at_2^{(11)}/\at_3^{(11)}/}
   \plot 52 33
   \plotfmla{\fourcrit\kappa_1/\nu/\at_3^{(11)}/\at_4^{(11)}/}
   \plot 52 23
   \plotfmla{\fourcrit\kappa_1/\nu/\at_4^{(11)}/\at_5^{(11)}/}
   \plot 52 13
   \plotfmla{\fourcrit\kappa_1/\nu/\at_5^{(11)}/\at_6^{(11)}/}
   \plot 52 3
   \plotfmla{\boxed{\kappa_2{\mapsto}\nu}}
   \plot 64 53
\plotfigend
\botcaption{Figure 3}
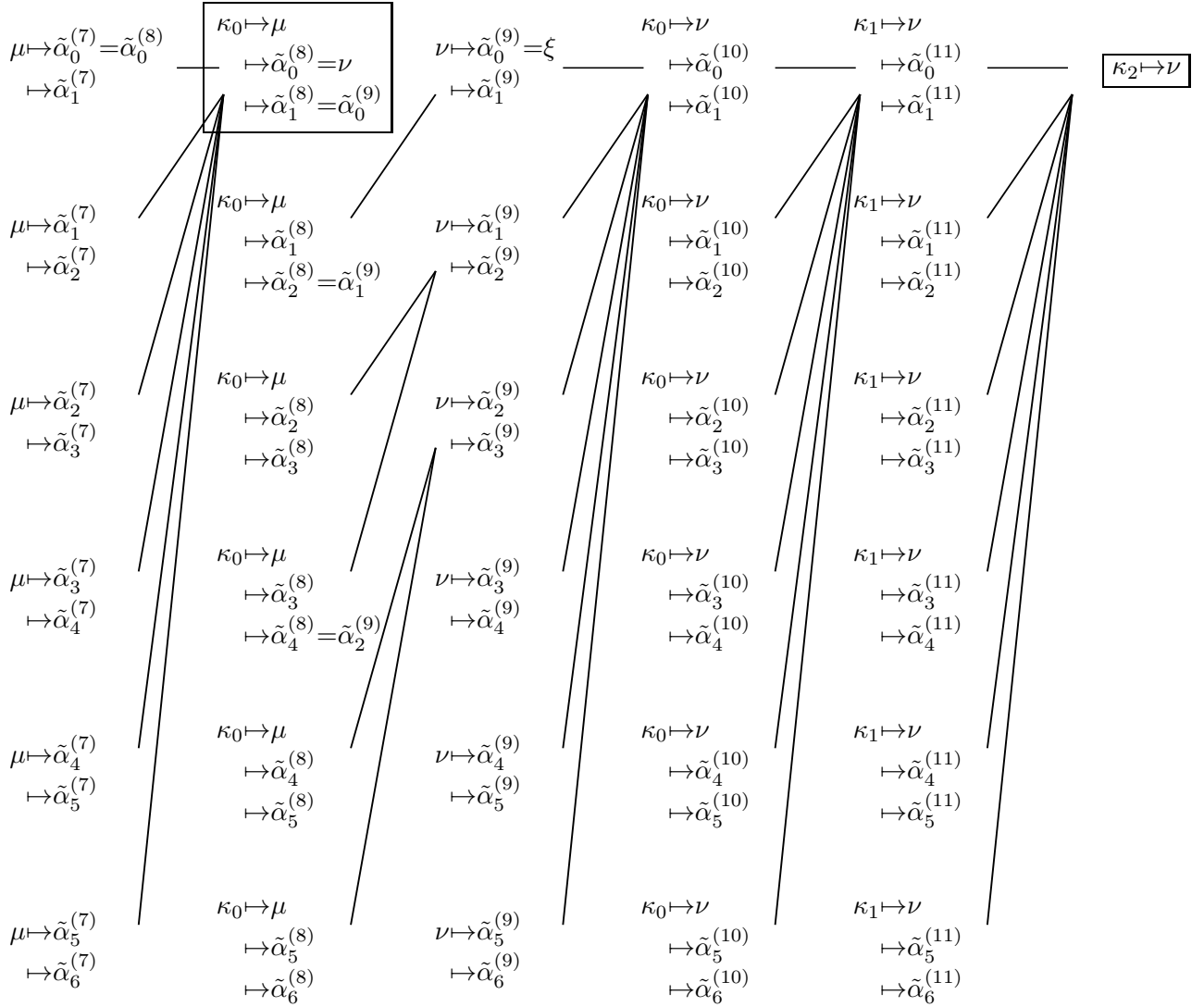
The embeddings $\et_m^{(n)}$ for $7 \le n \le 12$, $0 \le m \le 5$.
\endcaption
\endinsert

The new part of the construction is shown in Figure~3.  The boxed entries
are the given embeddings $k'$ and $k''$.  The new functions $\Ct_N$
indicating how many embeddings in column $N$ are obtained from a
given number of embeddings in column $N+1$ are defined as follows:
$$\gather \Ct_6(m) = 2^m; \\
\Ct_7(0)=0,\qquad \Ct_7(m+1)=\Ct_7(m) +
\Ctfunc17(\Ct_7(m))-8; \\
\Ct_8(m) = 2^m; \\
\allowdisplaybreak
\Ct_9(0)=0,\qquad \Ct_9(1)=\Ctfunc19(\Ctfunc16(1)-8)-8, \\
\Ct_9(m+1)=\Ct_9(m) +
\Ctfunc19(\Ct_9(m))-\Ctfunc16(1)\quad\text{for $m>0$;} \\
\Ct_{10}(0)=0,\qquad \Ct_{10}(m+1)=\Ct_{10}(m) +
\Ctfunc3{10}(\Ct_{10}(m)). \endgather$$

Instead of giving six more lemmas for the construction of columns
6--11, we will just give a description of how column number $N$ is
constructed from column $N+1$ for $N=6,7,8,9,10,11$.  From this,
it is straightforward to state corresponding lemmas, and the
proofs will be similar to those of preceding lemmas.

Column 6 is constructed from column 7 using Lemma~5;
no new results are needed.

If one has already used $n$ embeddings in column 8 to produce
$\Ct_7(n)$ embeddings in column 7 ending up with $\at^{(7)}_{\Ct_7(n)}
= \at^{(8)}_n$, then one can apply the preceding results to get
$\Ctfunc17(\Ct_7(n))$ embeddings in column~1, and one will have
$\at^{(1)}_8 = \mu$ and $\at^{(1)}_{\Ctfunc17(\Ct_7(n))} =
\at^{(8)}_n$.  Now apply $\et^{(8)}_n$ to all but the first $8$ of
these embeddings in column~1 to get $\Ctfunc17(\Ct_7(n))-8$ new
embeddings for column~7, for a total of $\Ct_7(n+1)$ embeddings ending
up with $\at^{(7)}_{\Ct_7(n+1)} = \at^{(8)}_{n+1}$.

Column 8 is constructed from column 9 using Lemma~5.

To build column 10 from column 9 requires an extra assumption not
visible in Figure~3: $\et^{(10)}_0(\mu) = \xi$.  Given this, one can
proceed as follows.  First, use the given embedding $k$ in column 6 to
produce $\Ctfunc16(1)$ embeddings in column~1, with $\at^{(1)}_8 = \mu$
and $\at^{(1)}_{\Ctfunc16(1)} = \nu$.  Apply $\et^{(10)}_0$ to all but
the first $8$ of these in order to get the first $\Ctfunc16(1)-8$
embeddings in column~9, with $\at^{(9)}_0 = \xi$ and
$\at^{(9)}_{\Ctfunc16(1)-8} = \at^{(10)}_0$.  These can be used to
produce $\Ctfunc19(\Ctfunc16(1)-8)$ embeddings in column~1, including
the ones produced before; apply $\et^{(10)}_0$ to all but
the first $8$ of these in order to get the first $\Ct_9(1)$
embeddings in column~9 (including
the ones produced before).  Now, given the $\Ct_9(n)$ embeddings
in column~9 produced from $n$ embeddings in column~10 ($n\ge1$),
do the previous columns' constructions to get
$\Ctfunc19(\Ct_9(n))$ embeddings in column 1, and apply
$\et^{(10)}_n$ to all but the first $\Ctfunc16(1)$ of these to get
new embeddings in column~9.

Building column 11 from column 10 requires the assumption
$\et^{(11)}_0(j)(\mu) = \xi$.  The construction is straightforward:
given $\Ct_{10}(n)$ embeddings in column 10, process the preceding
columns to get $\Ctfunc3{10}(\Ct_{10}(n))$ embeddings in column~3,
and apply $\et^{(11)}_n$ to all of these to get new embeddings
for column~10.

Finally, in building column 11, one does not have a true column 12
but just a single embedding~$k''$.  One gets column~11 by applying $k''$
to column~6: $\et^{(11)}_i = k''(\et^{(6)}_i)$ and
$\at^{(11)}_i = k''(\at^{(6)}_i)$ for all~$i$.  Since
$\Ctfunc6{11}(n) > n$ for all $n$, there is no problem continuing
this indefinitely to make column~11 as long as desired.

\topinsert
\centerline{\vbox{\offinterlineskip
\halign{
&\vrule#&\hfil$\,\vphantom{\boxed{\Ct_1\at_1^{(2)}}}#\,$\hfil\cr
\noalign{\hrule}
height2pt&\omit&&\omit&&\omit&&\omit&&\omit&&\omit&&\omit&&\omit&\cr
&n&&\text{Critical sequence of $\et_m^{(n)}$}&&\et_0^{(n)}&&\at_0^{(n)}&&
      \hat n&&\Ct_ n(0),(1)&&\Ct_ n(m{+}1)-\Ct_ n(m)&&
      \text{Growth}&\cr
height2pt&\omit&&\omit&&\omit&&\omit&&\omit&&\omit&&\omit&&\omit&\cr
\noalign{\hrule}
height2pt&\omit&&\omit&&\omit&&\omit&&\omit&&\omit&&\omit&&\omit&\cr
&0&&\at_m^{(0)} \mapsto \at_{m+1}^{(0)}&&j&&
      \kappa_0&&0&&1&&\Ct_ 0(m)&&F_3(m{-}3)&\cr
&1&&\kappa_0 \mapsto \at_m^{(1)} \mapsto \at_{m+1}^{(1)}&&j&&
      \kappa_1&&1&&1&&\Ct_ 1(m)&&F_3(m{-}3)&\cr
&2&&\kappa_1 \mapsto \at_m^{(2)} \mapsto \at_{m+1}^{(2)}&&&&
      \kappa_2&&1&&0&&\Ct_1(\Ct_2(m))&&F_4(m{-}4)&\cr
&3&&\kappa_0 \mapsto \kappa_1 \mapsto \at_m^{(3)} \mapsto \at_{m+1}^{(3)}&&j&&
      \kappa_2&&3&&1&&\Ct_3(m)&&F_3(m{-}3)&\cr
&4&&\kappa_2 \mapsto \at_m^{(4)} \mapsto \at_{m+1}^{(4)}&&&&
      \kappa_3&&1&&0,8&&\Ctfunc14(\Ct_4(m))-1&&F_5(m{-}2)&\cr
&5&&\kappa_0 \mapsto \kappa_2 \mapsto \at_m^{(5)} \mapsto \at_{m+1}^{(5)}&&&&
      \mu&&3&&0,2&&\Ctfunc35(\Ct_5(m))&&F_6(m{-}3)&\cr
&6&&\kappa_1 \mapsto \kappa_2 \mapsto \at_m^{(6)} \mapsto \at_{m+1}^{(6)}&&k&&
      \mu&&6&&1&&\Ct_6(m)&&F_3(m{-}3)&\cr
&7&&\mu \mapsto \at_m^{(7)} \mapsto \at_{m+1}^{(7)}&&&&
      \nu&&1&&0&&\Ctfunc17(\Ct_7(m))-8&&F_7(m{-}2)&\cr
&8&&\kappa_0 \mapsto \mu \mapsto \at_m^{(8)} \mapsto \at_{m+1}^{(8)}&&k'&&
      \nu&&8&&1&&\Ct_8(m)&&F_3(m{-}3)&\cr
&9&&\nu \mapsto \at_m^{(9)} \mapsto \at_{m+1}^{(9)}&&&&
      \xi&&1&&0,[{*}]&&\Ctfunc19(\Ct_9(m))-\Ctfunc16(1)&&F_8(m{-}2)&\cr
&10&&\kappa_0 \mapsto \nu \mapsto \at_m^{(10)} \mapsto \at_{m+1}^{(10)}&&&&
      k''(k)(\kappa_2)&&3&&0&&\Ctfunc3{10}(\Ct_{10}(m))&&F_9(m{-}2)&\cr
&11&&\kappa_1 \mapsto \nu \mapsto \at_m^{(11)} \mapsto \at_{m+1}^{(11)}&&&&
      k''(\mu)&&6&&0,\Ctfunc6{11}(?)&&&&&\cr
&12&&\kappa_2 \mapsto \nu&&k''&&
      &&&&1&&&&&\cr
height2pt&\omit&&\omit&&\omit&&\omit&&\omit&&\omit&&\omit&&\omit&\cr
\noalign{\hrule}
\noalign{\vskip 4pt}
\omit&\omit$\,[{*}] = \Ctfunc19(\Ctfunc16(1)-8)-8$\hidewidth\cr
}}}
\botcaption{Table 1}
Data for the construction of $\et_m^{(n)}$.
\endcaption
\endinsert

This entire construction is summarized
in Table~1, which includes: the beginning of the critical sequence
of $\et_m^{(n)}$; the initial boxed entry in column~$n$, if any;
the value of $\at_0^{(n)}$, usually used as an initial condition in the
construction; the number $\hat n$ such that column~$n$ is extended by
applying embeddings from column~$n{+}1$ to embeddings from column~$\hat n$;
the initial value for $\Ct_n$, or the initial two values if
necessary; the recursive formula for $\Ct_n(m+1)$ (which does not
apply to $m=0$ if two values were given in the preceding column); and
an estimate for the growth rate of
$\Ct_ n(m)$.

Straightforward inductions show that the values given under ``Growth''
are actually lower bounds on $\Ct_n(m)$.  In fact, we have
$\Ct_ n(m) = F_3(m-3)+3$ in the columns where $\Ct_n(m)=2^m$,
and $\Ct_ n(m) \ge F_i(m-d) + 4$ in the other columns, for the listed
values of $i$ and $d$; this trivial improvement makes the inductions
easier, and also makes it easier to compose the estimates in the
following calculations.

In order to conclude anything about $F(4)$ from this construction, we need
to know that certain critical points lie below $\kappa_4$.  For example,
once we know that $\xi < \kappa_4$, we can conclude that
$F(4) > \Ctfunc 08(1)$; the estimates in Table~1 allow us to
compute $\Ctfunc 08(1) = \Ctfunc07(\Ctfunc16(1)-\nobreak
8) > F_6(F_4(F_4(254)))$.
A better result can be obtained from the inequality
$k''(k''(\mu)) < \kappa_4$.  Since $k''(\mu) = k''(k(\kappa_2)) =
k''(k)(k''(\kappa_2)) = \et_0^{(11)}(\nu) = \at_0^{(11)} =
\at_{\Ctfunc36(1)}^{(10)}$, we get $\Ctfunc6{10}(\Ctfunc36(1))$ embeddings
in column~$6$ before reaching the ordinal $k''(\mu)$.  Therefore,
we get $\Ctfunc6{10}(\Ctfunc36(1))$ embeddings in column~$11$ before
reaching $k''(k''(\mu))$, so
$F(4) > \Ctfunc0{11}(\Ctfunc6{10}(\Ctfunc36(1))) > F_9(F_8(F_4(254)))$.

Finally, by pushing the computations in the next section a little farther,
we will show that $k''(k''(k)(k''(k)(k(j)(\mu)))) < \kappa_4$.  Since
$k(j)(\mu) = \at_1^{(5)}$, we get $\Ctfunc35(1)$ embeddings in column~$3$
before reaching $k(j)(\mu)$, and this leads to a total of
$\Ctfunc0{11}(\Ctfunc6{10}(\Ctfunc3{10}(\Ctfunc35(1))))$ critical points
produced below $k''(k''(k)(k''(k)(k(j)(\mu))))$.  This gives a lower
bound of $F_9(F_8(F_8(254)))$ for $F(4)$.  One could pursue
the inequalities further, but we have clearly reached the point
of diminishing returns.

There is no reason whatsoever to believe that the results in this
section are even close to optimal; the construction here is merely the
most successful from a series of {\it ad hoc} attempts to modify the
original construction to produce critical points below $\kappa_4$.  It
is probable that further modifications of the scheme would improve the
lower bound on $F(4)$.  One idea for modification would be to try to
produce a true column~$12$ instead of a single embedding, and it would
appear that the critical sequence for a member of this column should
begin with $\kappa_2 \mapsto \nu \mapsto \at_m^{(12)} \mapsto
\at_{m+1}^{(12)}$.
Unfortunately, this will not work with the present $k''$, because it
turns out that $k''(k''(\nu))$ is greater than $\kappa_4$.  So
apparently more drastic revisions would be needed.

\head 5. Critical point inequalities \endhead

In this section, we will use several basic methods to perform a number
of computations on critical points.  The main goal is to produce embeddings
$k$, $k'$, and $k''$ satisfying the assumptions stated in section~4, but
the accumulated facts may be useful for other purposes (e.g., as a place
to search for improved versions of the construction in section~4).

We will be concentrating on the sequence $\jsub n$ defined earlier by
$\jsub 1 = j$ and $\jsub {n+1} = \jsub n(j)$.  When looking at these
embeddings and their combinations, it will be helpful to omit some of
the parentheses in applications, thus writing $jj$ instead of $j(j)$.
In larger combinations, applications will be associated from the left:
$e_1e_2e_3$ means $(e_1e_2)e_3$, which is short for $(e_1(e_2))(e_3)$.
So one could write $jjjj$ for $\jsub 4$.

A great deal about the critical points of these embeddings can be
deduced from the simple rules $\crit{e_1e_2} = e_1(\crit{e_2})$ and
$e_1e_2(e_1(\beta)) = e_1(e_2(\beta))$.  For example, this sufficed in
section~1 to show that $\jsub 3(\kappa_1)$ is strictly between
$\kappa_2$ and $\kappa_3$.  (Since this ordinal will arise frequently,
it will be useful to have a name for it; let us call it
$\kappa_{2.5}$.)  Furthermore, all of the critical point manipulations
in sections~2 through~4 (except for the assumptions in section~4 which
were explicitly postponed until this section) used only these two
rules.

In order to proceed much further with the sequence $\jsub n$, though, we need
to use additional methods.  One very useful method is to explicitly use
the fact that $e_1e_2$ is obtained by applying~$e_1$ to initial segments
of $e_2$: $e_1e_2\restrict V_{e_1(\alpha)} = e_1(e_2 \restrict
V_\alpha)$.  Hence, if $e_2$ agrees with $e'_2$ up to~$\alpha$ (i.e.,
$e_2 \restrict V_\alpha = e'_2 \restrict V_\alpha$), then $e_1e_2$
agrees with $e_1e'_2$ up to $e_1(\alpha)$.  Furthermore, if $e_1$ agrees
with $e'_1$ up to some ordinal greater than the rank of
$e_2 \restrict V_\alpha$, then $e_1e_2$ agrees with $e'_1e'_2$
up to $e_1(\alpha)$.

Laver~\cite{5} has defined a variant form of `restriction' that turns out to
be more useful for these computations than the standard definition.  Let $e
\Lrestrict V_\beta = \{(x,y) \in V_\beta \times V_\beta \colon y \in e(x)\}$.
We again have $e_1e_2 \Lrestrict V_{e_1(\beta)} = e_1(e_2 \Lrestrict
V_\beta)$.  We can now define the modified version of `agreement up
to~$\beta$': say that $e \Lequiv\beta/ e'$ if $e \Lrestrict V_\beta = e'
\Lrestrict V_\beta$.  Then, for any limit ordinal $\beta$, $\Lequiv\beta/$
turns out to be an equivalence relation that respects composition and
application.  (To prove this, suppose that $e_1 \Lequiv\beta/ e'_1$ and $e_2
\Lequiv\beta/ e'_2$, and let $x,y \in V_\beta$; fix an ordinal $\gamma$
greater than the ranks of $x$ and $y$ but less than $\beta$.  If $y \in (e_1
\circ e_2)(x) = e_1(e_2(x))$, then $y \in e_1(z)$, where $z = e_2(x) \cap
V_\gamma$.  Since $e_1 \Lequiv\beta/ e'_1$, we have $y \in e'_1(z)$; since
$e_2 \Lequiv\beta/ e'_2$, we have $z = e'_2(x) \cap V_\gamma$.  Therefore, $y
\in (e'_1 \circ e'_2)(x)$.  For application, if $y \in e_1e_2(x)$, then $(x,y)
\in e_1(z)$, where $z = e_2 \Lrestrict V_\gamma$.  Now $e_1 \Lequiv\beta/
e'_1$ gives $(x,y) \in e'_1(z)$, and $e_2 \Lequiv\beta/ e'_2$ gives $z = e'_2
\Lrestrict V_\gamma$, so $y \in e'_1e'_2(x)$.)

The basic facts we will use about the relation $\Lequiv\beta/$ are:
\roster
\item $\Lequiv\beta/$ is an equivalence relation;
\item if $e_1 \Lequiv\beta/ e'_1$ and $e_2 \Lequiv\beta/ e'_2$, then
$e_1(e_2) \Lequiv\beta/ e'_1(e'_2)$ and $e_1\circ e_2 \Lequiv\beta/
e'_1\circ e'_2$;
\item if $e_2 \Lequiv\beta/ e'_2$, then
$e_1(e_2) \Lequiv e_1(\beta)/ e_1(e'_2)$;
\item $e \Lequiv\crit e/ \id$, where $\id$ is the identity function; and
\item if $e
\Lequiv\beta/ e'$ and either $e(\alpha)$ or $e'(\alpha)$ is less than $\beta$,
then $e(\alpha) = e'(\alpha)$.
\endroster
The first two of these were proved above, and the third is a consequence
of the elementarity of $e_1$; the last two are easy from the
definitions.  [For \therosteritem5, we cannot
have, for example, $e(\alpha)<\beta$ and $e(\alpha) < e'(\alpha)$,
because then the pair $(\alpha,e(\alpha))$ would be in
$e' \Lrestrict V_\beta$ but not in $e \Lrestrict V_\beta$.]

One application of these restriction methods is the following lemma.

\proclaim{Lemma 15} If\/ $\crit e \ge \kappa_2$, then $ejjj(\kappa_1) =
e(\kappa_{2.5})$. \endproclaim

\demo{Proof 1} Note that, if $e(x) = x$, then $ee'(x) = ee'(e(x)) =
e(e'(x))$.  Now, since $\kappa_2$ is regular, the map $j \restrict
\kappa_1$ from $\kappa_1$ to $\kappa_2$ is not cofinal in $\kappa_2$, so
the rank of $j \restrict V_{\kappa_1}$ is less than $\kappa_2$. Since
$ej(\kappa_1) = e(j(\kappa_1)) = e(\kappa_2) \ge \kappa_2$, $\kappa_1$
and $j\restrict \kappa_1$ are in the domain of $ej(j \restrict
V_{\kappa_1})$, and we have $ej(j \restrict V_{\kappa_1})(\kappa_1) =
e(j(j \restrict V_{\kappa_1})(\kappa_1)) = e(\kappa_2) > \kappa_1$.
Therefore, $$ejjj(\kappa_1) = ej(j \restrict V_{\kappa_1})(j \restrict
V_{\kappa_1})(\kappa_1) = e(j(j \restrict V_{\kappa_1})(j \restrict
V_{\kappa_1})(\kappa_1)) = e(jjj(\kappa_1)) = e(\kappa_{2.5}).$$
\enddemo

\demo{Proof 2} Since $\crit e \ge \kappa_2$, we have $e
\Lequiv\kappa_2/ \id$.  This gives $ej \Lequiv\kappa_2/ \id(j) = j$ and
hence $ejj \Lequiv\kappa_2/ jj$, so $$ej(\kappa_2) = ej(j(\kappa_1)) =
ejj(ej(\kappa_1)) \ge ejj(\kappa_2) = ejjj(ejj(\kappa_1)) \ge
ejjj(\kappa_2) > ejjj(\kappa_1).$$  We now have $ejjj \Lequiv
ejj(\kappa_2)/ ejj(ej) \Lequiv ej(\kappa_2)/ ej(ej)(ej) = e(jjj)$.
Therefore, $ejjj(\kappa_1) = e(jjj)(\kappa_1) = e(jjj)(e(\kappa_1)) =
e(jjj(\kappa_1)) = e(\kappa_{2.5})$. \QED

From now on, the derivations in this section will be in the form of
Proof~2 (using modified restrictions) rather than that of Proof~1
(using standard restrictions).  Derivations using standard restrictions
are sometimes easier to come up with---the above lemma is an example---but
derivations using modified restrictions give a little more information
and are applicable in other contexts besides that of elementary embeddings
(as in forthcoming work of Jech and myself~\cite{3}).

Using this lemma and the basic rules, we can now derive more information
about the embeddings~$\jsub n$ for relatively small $n$, as shown in
Table~2.  In this table and from now on, $\kappa_i^n$ is an abbreviation
for $\jsub n(\kappa_i)$.

\midinsert
\centerline{
\hfill
\vbox{\offinterlineskip
\def\\#1//{\rlap{$#1$}\hphantom{\kappa_{2.5}}}
\halign{
\vrule#&\hfil$\,\vphantom{\boxed{\kappa_2^7}}#\,$\hfil&
      \vrule#&$\,#\,$\hfil&\vrule#&$\quad#\quad$\hfil&\vrule#\cr
\noalign{\hrule}
height2pt&\omit&&\omit&&\omit&\cr
&n&&\omit\hfil Critical sequence of $\jsub n$\hfil&&
      \omit\hfil\,Other values\,\hfil&\cr
height2pt&\omit&&\omit&&\omit&\cr
\noalign{\hrule}
height2pt&\omit&&\omit&&\omit&\cr
&1&&\\\kappa_0// \mapsto \\\kappa_1// \mapsto \\\kappa_2// \mapsto
      \\\kappa_3// \mapsto \kappa_4&&&\cr
&2&&\\\kappa_1// \mapsto \\\kappa_2// \mapsto \\\kappa_3// \mapsto 
      \kappa_4&&&\cr
&3&&\\\kappa_0// \mapsto \\\kappa_2// \mapsto \\\kappa_3// \mapsto \kappa_4&&
      \kappa_1 \mapsto \kappa_{2.5}&\cr
&4&&\\\kappa_2// \mapsto \\\kappa_{2.5}// \mapsto \\\kappa_3// \mapsto
      \kappa_4&&&\cr
&5&&\\\kappa_0// \mapsto \\\kappa_1// \mapsto \\\kappa_{2.5}// \mapsto
      \kappa_4&&&\cr
&6&&\\\kappa_1// \mapsto \\\kappa_{2.5}// \mapsto \\\kappa_2^5// \mapsto
      {>}\kappa_4&&&\cr
&7&&\\\kappa_0// \mapsto \\\kappa_{2.5}// \mapsto \\\kappa_2^6// \mapsto 
      \\\kappa_3^6// \mapsto {>}\kappa_4&&\kappa_1 \mapsto \kappa_3&\cr
&8&&\\\kappa_{2.5}// \mapsto \\\kappa_3// \mapsto \\\kappa_2^7// \mapsto
      \\\kappa_3^7// \mapsto {>}\kappa_4&&&\cr
height2pt&\omit&&\omit&&\omit&\cr
\noalign{\hrule}
}}
\hfill
\vbox{\offinterlineskip
\def\\#1//{\rlap{$#1$}\hphantom{\kappa_{2.5}}}
\halign{
\vrule#&\hfil$\,\vphantom{\boxed{\kappa_2^7}}#\,$\hfil&
      \vrule#&$\,#\,$\hfil&\vrule#&$\quad#\quad$\hfil&\vrule#\cr
\noalign{\hrule}
height2pt&\omit&&\omit&&\omit&\cr
&n&&\omit\hfil Critical sequence of $\jsub n$\hfil&&
      \omit\hfil\,Other values\,\hfil&\cr
height2pt&\omit&&\omit&&\omit&\cr
\noalign{\hrule}
height2pt&\omit&&\omit&&\omit&\cr
&9&&\\\kappa_0// \mapsto \\\kappa_1// \mapsto \\\kappa_2// \mapsto
      \\\kappa_2^7// \mapsto {>}\kappa_4&&&\cr
&10&&\\\kappa_1// \mapsto \\\kappa_2// \mapsto \\\kappa_2^7// \mapsto
      \\\kappa_3^9// \mapsto {>}\kappa_4&&&\cr
&11&&\\\kappa_0// \mapsto \\\kappa_2// \mapsto \\\kappa_2^7// \mapsto
      \\\kappa_3^{10}// \mapsto {>}\kappa_4&&\kappa_1 \mapsto \kappa_3&\cr
&12&&\\\kappa_2// \mapsto \\\kappa_3// \mapsto \\\kappa_2^7// \mapsto
      \\\kappa_3^{11}// \mapsto {>}\kappa_4&&&\cr
&13&&\\\kappa_0// \mapsto \\\kappa_1// \mapsto \\\kappa_3// \mapsto
      \\\kappa_2^7// \mapsto {>}\kappa_4&&&\cr
&14&&\\\kappa_1// \mapsto \\\kappa_3// \mapsto \\\kappa_2^{13}// \mapsto
      \\\kappa_2^7// \mapsto {>}\kappa_4&&&\cr
&15&&\\\kappa_0// \mapsto \\\kappa_3// \mapsto \\\kappa_2^{14}// \mapsto
      \\\kappa_2^{13}// \mapsto {>}\kappa_4&&&\cr
&16&&\\\kappa_3// \mapsto \\\kappa_1^{15}// \mapsto \\\kappa_2^{15}// \mapsto
      \\\kappa_2^{14}// \mapsto {>}\kappa_4&&&\cr
height2pt&\omit&&\omit&&\omit&\cr
\noalign{\hrule}
}}
\hfill}
\botcaption{Table 2}
Selected values of $\jsub n$ for $1 \le n \le 16$.
\endcaption
\endinsert

The `${>}\kappa_4$' entry for $\jsub 6$ comes from $\kappa_3^5 >
\kappa_{2.5}^5 = \kappa_4$.  For the `${>}\kappa_4$' entry for $\jsub
n$ ($7 \le n \le 16$), we need the fact that $\jsub {n-1}(\kappa_4) >
\kappa_4$.  An easy way to see this is to prove by induction that an
ordinal fixed by an embedding $e$ which lies above the critical point
of $e$ must lie above the entire critical sequence of $e$.  Another way
is to conclude from a theorem of Kunen (in the form given in
Solovay-Reinhardt-Kanamori~\cite{8}) or the results of Laver and Steel
mentioned previously that the critical sequence of $e$ must be cofinal
in $\lambda$.  (Recall that $j\colon V_\lambda \to V_\lambda$.)

Since Laver~\cite{6} has shown that there are exactly $m$ critical
points below $\crit{\jsub {2^m}}$ for all $m$, the fact that $\kappa_3
= \crit{\jsub {16}}$ implies that there are no critical points below
$\kappa_3$ other than the four we already know about.

We now have enough information to prove a number of the facts needed
for section~4.  Let $k = \jsub {10}$; then \fourseq k/
\kappa_1 / \kappa_2 / \mu / \nu/, where $\mu =
\kappa_2^7$ and $\nu = \kappa_3^9$.  Also, $k(j)(\kappa_1) = \jsub
{11}(\kappa_1) = \kappa_3$.  Next, let $k' = \jsub {10}(\jsub {11})$;
\fourseq k'/\jsub {10}(\kappa_0) / \jsub
{10}(\kappa_2) / \jsub {10}(\kappa_2^7) / \jsub
{10}(\kappa_3^{10})/, which is $\kappa_0 \mapsto \mu \mapsto \nu \mapsto \xi$,
where $\xi = \jsub {10}(\kappa_3^{10})$.  Finally, if we let $k'' =
\jsub 9(\jsub {14})$, then the critical sequence of $k''$ will begin with
$\jsub 9(\kappa_1) \mapsto \jsub 9(\kappa_3)$, which is $\kappa_2 \mapsto
\nu$.  All that remains is to show that $k''(k)(j)(\mu) = \xi$ and that
the relevant critical points lie below $\kappa_4$.

The main problem we will have to solve is that of comparing various
critical points.  Certain such comparisons can be determined directly
from Table~2; for example, $\kappa_2^7 < \kappa_2^6$ because
$\kappa_2^6 = \kappa_{2.5}^7$.  A thorough examination of Table~2 gives
the following list of critical points of the forms $\kappa_i$ and $\kappa_i^n$:
$$\multline
\kappa_0 < \kappa_1 < \kappa_2 < \kappa_{2.5} < \kappa_3 < \kappa_1^{15}
< \kappa_2^{15} < \kappa_{2.5}^{15} < \kappa_2^{14} < \kappa_{2.5}^{14}
< \kappa_2^{13} < \kappa_{2.5}^{13} \\ < \kappa_2^7 < \kappa_2^6 <
\kappa_2^5 < \kappa_4 < \kappa_3^5.
\endmultline$$
However, the ordinals $\kappa_{2.5}^n$ for $n = 9,10,11$ and $\kappa_3^n$
for $n = 6,7,9,10,11$ cannot yet be placed in this list; we only have
the partial information that $\kappa_2^7 < \kappa_{2.5}^{11} <
\kappa_3^{11} < \kappa_3^{10}$, $\kappa_2^7 < \kappa_{2.5}^{10} <
\kappa_3^{10} < \kappa_3^9$, $\kappa_2^7 < \kappa_{2.5}^9 < \kappa_3^9$,
$\kappa_2^6 < \kappa_3^7 < \kappa_3^6$, and $\kappa_2^5 < \kappa_3^6$.
The critical points $\kappa_{2.5}^n$ for $n = 1,2,3,12$ are not mentioned
here because, as we will see soon, they already occur in the above list.
We will return later to the problem of sorting the remaining critical points
into the correct order.

One useful method for critical point computations is to approximate
(in the sense of $\Lequiv\alpha/$) complicated expressions such as
$\jsub n$ for relatively large $n$ by simpler expressions.  In order to do
this, some basic computations are useful.

First, from Laver~\cite{5},
for any embeddings $e$ and $e_1,\dots,e_l$, we have
$ee_1e_2\dots e_l \Lequiv\theta/ e(e_1e_2\dots e_l)$, where
$\theta$ is the minimum of $ee_1e_2\dots e_i(\crit e)$ for
$1 \le i < l$.  This is proved by induction:
$$ee_1e_2\dots e_l \Lequiv ee_1e_2\dots e_{l-1}(\crit e)/ ee_1e_2\dots e_{l-1}
(ee_l) \Lequiv\theta/ e(e_1e_2\dots e_{l-1})(ee_l) = e(e_1e_2\dots e_l).$$

Second, we have $e \Lequiv e(\crit{e'})/ e \circ e'$.  To see this,
note that $\crit{ee'} = e(\crit{e'})$, so $\id \Lequiv e(\crit{e'})/
ee'$, so $e = \id \circ e \Lequiv e(\crit{e'})/ ee' \circ e = e \circ
e'$.  A special case of this is $ee' \Lequiv ee'(\crit e)/ ee' \circ e
= e \circ e'$.  This is a version of the fact given earlier, that $e(x)
= x$ implies $ee'(x) = ee'(e(x)) = e(e'(x))$.

Using these facts, we can see that $\jsub 3$ is approximated by $j
\circ j$, because $\jsub 3 \Lequiv\kappa_{2.5}/ \jsub 2 \circ j = j
\circ j$.  But a glance at Table~2 indicates that $j\circ j$ should
approximate $\jsub {11}$ better than it approximates~$\jsub 3$, since
$\jsub {11}(\kappa_1) = (j \circ j)(\kappa_1) = \kappa_3 \ne \jsub
3(\kappa_1)$.  To see that this is indeed the case, note that $\jsub
{11} = \jsub 8jjj \Lequiv\kappa_2^7/ \jsub 8(jjj) = \jsub 8\jsub 3$,
since both $\jsub 8j(\kappa_{2.5})$ and $\jsub 8jj(\kappa_{2.5})$ are
greater than $\kappa_2^7$.  Similarly, we get $\jsub 8
\Lequiv\kappa_2^7/ \jsub 4\jsub 4 = \jsub 3\jsub 2$, so $\jsub {11}
\Lequiv\kappa_2^7/ \jsub 3\jsub 2\jsub 3$.  Now note that $$\jsub
3\jsub 2\jsub 3(\kappa_{2.5}) > \jsub 3\jsub 2\jsub 3(\kappa_2) = \jsub
3\jsub 2(\jsub 3(\kappa_2)) = \jsub 3(\kappa_3) = \kappa_4 >
\kappa_2^7,$$ so $\jsub {11} \Lequiv\kappa_2^7/ \jsub 3\jsub 2\jsub 3
\Lequiv\kappa_2^7/ \jsub 3\jsub 2 \circ \jsub 3 = \jsub 3 \circ \jsub 2
= \jsub 2 \circ j = j \circ j$.  Since $\jsub {11}(\kappa_2) =
\kappa_2^7$, we get as a consequence that $\jsub {11} \restrict
V_{\kappa_2} = (j \circ j) \restrict V_{\kappa_2}$.

Similarly, one can approximate $\jsub {12}$ by $j\jsub 2$, since $\jsub
{12} = \jsub {11}j \Lequiv\kappa_2^7/ (j \circ j)j = j\jsub 2$.
Continuing, we get the following:  $$\gather \jsub {13} = \jsub {12}j
\Lequiv\kappa_2^7/ j\jsub 2j \\ \jsub {14} = \jsub {12}jj
\Lequiv\kappa_2^{13}/ \jsub {12}(jj) \Lequiv\kappa_2^7/ j\jsub 2(jj) =
j\jsub 3 \\ \jsub {15} = \jsub {12}jjj \Lequiv\kappa_2^{14}/ \jsub
{12}(jjj) \Lequiv\kappa_2^7/ j\jsub 2(jjj) = \jsub 2\jsub 2(\jsub 2j) =
\jsub 2\jsub 3.\endgather$$

These approximations allow us to evaluate certain critical points in
terms of others.  First, since $\jsub {14}(\kappa_2) < \kappa_2^{13}$,
we have $\kappa_2^{14} = j\jsub 3(\kappa_2) = j\jsub 3(j(\kappa_1)) =
j(\jsub 3(\kappa_1)) = \kappa_{2.5}^1$.  Similarly, $\kappa_2^{15} <
\kappa_2^{14}$ implies $\kappa_2^{15} = \jsub 2\jsub 3(\kappa_2) =
\jsub 2\jsub 3(\jsub 2(\kappa_1)) = \kappa_{2.5}^2$.  Since
$\kappa_1^{15} < \kappa_2^{14}$, we have $\kappa_1^{15} = \jsub
{12}(jjj)(\kappa_1) = \jsub {12}\jsub 3(\jsub {12}(\kappa_1)) =
\kappa_{2.5}^{12}$.  But we can go further here; since $\jsub {12}
\Lequiv\kappa_2^7/ j\jsub 2 = \jsub 2\jsub 2 = \jsub 3\jsub 3$, we get
$\kappa_{2.5}^{12} = \jsub 3\jsub 3(\kappa_{2.5}) = \jsub 3\jsub
3(\jsub 3(\kappa_1)) = \kappa_{2.5}^3$.  So the critical points
$\kappa_{2.5}^n$ for $n=1,2,3,12$ are already in our list and do not
need to be added.

We can now get the approximations $\jsub {13} \Lequiv\kappa_2^{13}/ j
\circ \jsub 2$ and $\jsub {14} \Lequiv\kappa_2^{14}/ \jsub 2 \circ
\jsub 2$ as follows:  $$\gather \jsub {13} \Lequiv \jsub
{13}(\crit{j\jsub 2})/ \jsub {13} \circ j\jsub 2 \Lequiv\kappa_2^7/
j\jsub 2j \circ j\jsub 2 = j\jsub 2 \circ j = j \circ \jsub 2, \\ \jsub
{14} \Lequiv\kappa_2^{13}/ j\jsub 3 \Lequiv j(\kappa_{2.5})/ j(j \circ
j) = \jsub 2 \circ \jsub 2. \endgather$$ The latter implies that $\jsub
{15} \Lequiv\kappa_1^{15}/ j\circ j\circ j$, since $$\jsub {15} = \jsub
{14}j \Lequiv\kappa_1^{15}/ \jsub {14} \circ j \Lequiv\kappa_2^{14}/
\jsub 2 \circ \jsub 2 \circ j = \jsub 2 \circ j \circ j = j \circ j
\circ j.$$  This is not as good an approximation for $\jsub {15}$ as
$\jsub 2\jsub 3$, but it will suffice for the following.

We are now ready to show that $k''(k)(j)(\mu) = \xi$.  Since $k'' =
\jsub 9\jsub {14}$ and $k = \jsub {10} = \jsub 9j$, we have $k''k =
\jsub 9\jsub {15}$.  We also have $\jsub 9\jsub {15} \Lequiv \jsub
9(\kappa_1^{15})/ \jsub 9(j \circ j \circ j)$; since $$\multline \jsub
9(\jsub {15}(\kappa_1)) = \jsub 9\jsub {15}(\jsub 9(\kappa_1)) = \jsub
9\jsub {15}(\kappa_2) = \jsub 9\jsub {15}(j(\kappa_1)) = \jsub 9\jsub
{15}j(\jsub 9\jsub {15}(\kappa_1)) \\ = \jsub 9\jsub {15}j(\jsub 9\jsub
{15}(\jsub 9(\kappa_0))) = \jsub 9\jsub {15}j(\kappa_3^9) > \jsub
9\jsub {15}j(\mu),\endmultline$$ we get $$\multline \jsub 9\jsub {15}j(\mu) =
\jsub 9(j\circ j\circ j)j(\mu) = (\jsub {10}\circ \jsub {10}\circ \jsub
{10})j(\mu) \\ = \jsub {10}(\jsub {10}(\jsub {11}))(\jsub {10}(\jsub
{10}(\kappa_1))) = \jsub {10}(\jsub {10}(\kappa_3)) = \xi. \endmultline$$

It remains to show that a certain critical point lies below $\kappa_4$.
In order to do so, we will study a more general question: Given an
embedding $e$ and a critical point $\gamma$, what is the least $\alpha$
such that $e(\alpha) \ge \gamma$?  (This will be solved here only in a
few special cases; in general it seems quite difficult.)  This ordinal
$\alpha$ is often not a critical point, so we will start by studying
some additional ordinals which are definable from the embeddings in
$\jalgp$.

For any embedding $e$ and ordinal $\alpha$, let $e({<}\alpha)$ be the
strict supremum of $e(\beta)$ for $\beta < \alpha$ (i.e., the least
ordinal greater than all such $e(\beta)$).  We clearly have
$e({<}\alpha) \le e(\alpha)$.  If $\rho$ is the cofinality of $\alpha$,
then $e(\rho) > \rho$ implies $e({<}\alpha) < e(\alpha)$, because the
cofinality of $e({<}\alpha)$ is $\rho$ while that of $e(\alpha)$ is
$e(\rho)$.  (Conversely, if $e(\rho) = \rho$, then $e({<}\alpha) =
e(\alpha)$.  To see this, let $g\colon \rho \to \alpha$ be increasing
and cofinal in $\alpha$; then $e(g)$ must be increasing and cofinal in
$e(\alpha)$.  But the domain of $e(g)$ is $e(\rho) = \rho$, and
$e(g)(e(\delta)) = e(g(\delta))$ for $\delta<\rho$, so the ordinals
$e(\beta)$ for $\beta<\alpha$ are cofinal in $e(\alpha)$.)

As usual, elementarity implies that application distributes over this
new operation: $e'(e({<}\alpha)) = e'e({<}e'(\alpha))$.  Also, for any
$e'$, $e$, and $\alpha$, the ordinals $e(\beta)$ for $\beta<\alpha$ are
cofinal in $e({<}\alpha)$, so the ordinals $e'(e(\beta))$ for
$\beta<\alpha$ are cofinal in $e'({<}e({<}\alpha))$; this gives the
formula $e'({<}e({<}\alpha)) = (e'\circ e)({<}\alpha)$.
Another useful fact is that $e \Lequiv e(\alpha)/ e'$ implies $e(\beta)
= e'(\beta)$ for $\beta<\alpha$, and hence $e({<}\alpha) =
e'({<}\alpha)$.

Two of these ordinals arise frequently enough that it is useful to have
short names for them: let $\sigma_1 = j({<}\kappa_1)$  and $\sigma_2 =
j({<}\kappa_2)$.  Then we can see from the above that $\sigma_1$ is
strictly between $\kappa_1$ and~$\kappa_2$, while $\sigma_2$ is strictly
between $\kappa_2$ and $\kappa_3$.  To see where $\sigma_2$ lies
relative to $\kappa_{2.5}$, use the following computation:
$$\sigma_2 = j({<}\kappa_2) > j(\sigma_1) = \jsub 2({<}\kappa_2) >
\jsub 2(\sigma_1) = \jsub 3({<}\kappa_2) > \jsub 3(\kappa_1) = \kappa_{2.5}.$$
On the other hand, the fact that $\jsub 3 \Lequiv\kappa_{2.5}/ j\circ j$
yields $$j({<}\sigma_1) = (j\circ j)({<}\kappa_1) = \jsub 3({<}\kappa_1) <
\jsub 3(\kappa_1) = \kappa_{2.5}.$$  So $\kappa_{2.5}$ lies within a gap in
the range of $j$, and $\sigma_1$ is the least ordinal that $j$ sends
above $\kappa_{2.5}$.

The inequality $(j\circ j)({<}\kappa_1) < \kappa_{2.5}$ will be very
useful later, because we can apply an arbitrary embedding $e$ to
it to get $(ej \circ ej)({<}e(\kappa_1)) < e(\kappa_{2.5})$.  In other
words, $$\text{if}\quad \beta < e(\kappa_1),\quad\text{then}\quad
ej(ej(\beta)) < e(\kappa_{2.5}). \tag\dag$$

We can now try to determine the ordinals at which various embeddings
first jump beyond $\kappa_4$.  Let us first consider $\jsub 6$; the
ordinal here turns out to be $\jsub 5(\sigma_1)$.  To see this, compute
as follows:  $$\gather \jsub 6(\jsub 5(\sigma_1)) = \jsub
5(j(\sigma_1)) > \jsub 5(\kappa_{2.5}) = \kappa_4, \\ \jsub 6({<}\jsub
5(\sigma_1)) = \jsub 5(j({<}\sigma_1)) < \jsub 5(\kappa_{2.5}) =
\kappa_4.\endgather$$ In particular, since $\jsub 5(\sigma_1) = \jsub
6({<}\kappa_{2.5}) > \kappa_2^6$, we get $\kappa_3^6 < \jsub
6(\kappa_2^6) < \kappa_4$.

[One could give a proof of $\jsub 6(\kappa_2^6) < \kappa_4$ without
using ordinals other than critical points, as follows: we have $\jsub 3
\Lequiv\kappa_{2.5}/ j\circ j$, so $\jsub 5\jsub 3 \Lequiv\kappa_4/
\jsub 5(j\circ j)$, so $$\kappa_4 = \jsub 5(\jsub 3(\kappa_1)) = \jsub
5\jsub 3(\kappa_{2.5}) > \jsub 5\jsub 3(\kappa_2) = \jsub 5(j \circ
j)(\kappa_2) = \jsub 6(\jsub 6(\kappa_2)).$$ In fact, Jech and I have
shown~\cite{3} that, in theory, any argument involving the ordinals
$e({<}\alpha)$ and the methods used here can be translated into an
argument referring only to critical points.  However, such arguments
tend to involve long, roundabout, and opaque computations.]

The embedding $\jsub 9$ is a more difficult case.  The least ordinal
sent above $\kappa_4$ by $\jsub 9$ is $\jsub 8(\sigma_2)$.  To see that
$\jsub 9(\jsub 8(\sigma_2)) > \kappa_4$, note that $$\multline
j(\sigma_2) = \jsub 2({<}\kappa_3) > \jsub 2(\sigma_2) = \jsub
3({<}\kappa_3) > \jsub 3(\sigma_2) = \jsub 4({<}\kappa_3) > \jsub
4(\sigma_2) = \jsub 5({<}\kappa_{2.5}) \\ > \jsub 5(j({<}\sigma_1)) =
\jsub 6({<}\jsub 5(\sigma_1)) > \kappa_3^6 > \kappa_3^7,\endmultline$$
so $\jsub 9(\jsub 8(\sigma_2)) = \jsub 8(j(\sigma_2)) > \jsub
8(\kappa_3^7) > \kappa_4$.  For the other side, first use the fact that
$\jsub {10} \Lequiv\kappa_2^7/ \jsub 8(\jsub 2)$ (since $\kappa_{2.5}^9
> \kappa_2^9 = \kappa_2^7$) to see that $\jsub {10}(\sigma_1) = \jsub
8(jj)(\sigma_1) = \jsub 8(\jsub 2(\sigma_1))$ (since $\sigma_1 <
\kappa_2$ and $\kappa_2^{10} = \kappa_2^7$).  From this, we get
$$\multline \jsub 9({<}\jsub 8(\sigma_2)) = \jsub 8(j({<}\sigma_2)) =
\jsub 8((j\circ j)({<}\kappa_2)) = \jsub 8(\jsub {11}({<}\kappa_2))
\\ = \jsub 8(\jsub {10}(\sigma_1)) = \jsub 8(\jsub 8(\jsub
2(\sigma_1))). \endmultline$$ Since $\jsub 2(\alpha_1) < \kappa_3 =
\jsub 7(\kappa_1)$, we can now apply \thetag{\dag} to get $\jsub
8(\jsub 8(\jsub 2(\sigma_1))) < \jsub 7(\kappa_{2.5}) = \kappa_2^6 <
\kappa_4$.

By methods similar to the above (with each $n$ handled separately), it
can be shown that: for $7 \le n \le 13$, the least $\alpha$ such that
$\jsub n(\alpha) > \kappa_4$ is $\jsub {n-1}(\sigma_2)$; for $n=14$ or $n=15$,
the least such $\alpha$ is $\jsub {n-1}(j(\sigma_1))$.

We now have the needed tools to complete the ordering of the critical
points $\kappa_i^n$ for $i \le 3$ and $n \le 16$.  The main tool needed
here is the formula $j({<}\sigma_1) < \kappa_{2.5}$; as with
\thetag{\dag}, we can apply any embedding $e$ to this to get
$ej({<}e(\sigma_1)) < e(\kappa_{2.5})$.  For instance, take $e = \jsub
9$ to get $\jsub {10}({<}\jsub 9(\sigma_1)) < \kappa_{2.5}^9$; from
this and the calculation $$\jsub 9(\sigma_1) = \jsub {10}({<}\kappa_2)
> \jsub {10}(\sigma_1) = \jsub {11}({<}\kappa_2) > \jsub {11}(\sigma_1)
= \jsub {12}({<}\kappa_3) > \kappa_2^{12} = \kappa_3,$$ we get
$\kappa_3^{10} < \kappa_{2.5}^9$.  The same argument shows that
$\kappa_3^{11} < \kappa_{2.5}^{10}$; it also shows that $\kappa_3^7 <
\kappa_{2.5}^6 = \kappa_2^5$, once we verify that $\jsub 6(\sigma_1) >
\kappa_3$.  To see this, note that $\jsub 4jj \Lequiv\kappa_2^5/ \jsub
4(jj)$ and $\kappa_2^5 > \kappa_2^6 > \jsub 6(\sigma_1)$, so $\jsub
6(\sigma_1) = \jsub 4\jsub 2(\sigma_1) = \jsub 4(\jsub 2(\sigma_1)) >
\jsub 4(\kappa_{2.5}) = \kappa_3$.

It remains to compare $\kappa_3^9$ with $\kappa_2^6$; we can show that
$\kappa_3^9 < \kappa_2^6$ as follows.  Apply $\jsub 7$ to the formula
$j({<}\sigma_1) < \kappa_{2.5}$ to get $\jsub 8({<}\jsub 7(\sigma_1)) <
\kappa_{2.5}^7 = \kappa_2^6$.  Since $\jsub 9(\kappa_3) = \jsub 9(\jsub
8(\kappa_{2.5})) = \jsub 8(j(\kappa_{2.5}))$, all we need to show is
that $j(\kappa_{2.5}) < \jsub 7(\sigma_1)$.  This is proved by the
following computation, part of which was already done in the preceding
paragraph:  $$\multline \jsub 7(\sigma_1) = \jsub 8({<}\kappa_3) >
\jsub 8(\sigma_2) = \jsub 9({<}\kappa_2) > \jsub 9(\sigma_1) \\ > \jsub
{12}({<}\kappa_3) > \jsub {12}(\sigma_2) = \jsub {13}({<}\kappa_3) >
\kappa_{2.5}^{13} > \kappa_2^{14} = \kappa_{2.5}^1.\endmultline$$ This
completes the comparison of these critical points $\kappa_i^n$; the
part of the list from $\kappa_2^7$ on is $$\kappa_2^7 <
\kappa_{2.5}^{11} < \kappa_3^{11} < \kappa_{2.5}^{10} < \kappa_3^{10} <
\kappa_{2.5}^9 < \kappa_3^9 < \kappa_2^6 < \kappa_3^7 < \kappa_2^5 <
\kappa_3^6 < \kappa_4 < \kappa_3^5.$$

We are finally ready to prove the last assumption in section~4,
$k''(k''(k)(k''(k)(k(j)(\mu)))) < \kappa_4$.  We have $k'' = \jsub
9\jsub {14}$, $k''(k) = \jsub 9\jsub {15}$, and $k(j) = \jsub {11}$.
This gives $k(j)(\mu) = \jsub {11}(\jsub {10}(\kappa_2)) = \jsub
{10}(j(\kappa_2)) = \kappa_3^{10}$.  We saw above that $\kappa_3^{10} <
\kappa_{2.5}^9$.  Since $\jsub {14}(\kappa_1) = \kappa_3 >
\kappa_{2.5}$, \thetag{\dag} implies that $\jsub {15}(\jsub
{15}(\kappa_{2.5})) < \jsub {14}(\kappa_{2.5})$; we can apply $\jsub 9$
to this to get $\jsub 9\jsub {15}(\jsub 9\jsub {15}(\kappa_{2.5}^9)) <
\jsub 9(\kappa_{2.5}^{14})$.  In the same way, since $\jsub
{13}(\kappa_1) > \kappa_{2.5}$, \thetag{\dag} yields $\jsub {14}(\jsub
{14}(\kappa_{2.5})) < \jsub {13}(\kappa_{2.5})$ and hence $\jsub 9\jsub
{14}(\jsub 9(\kappa_{2.5}^{14})) < \jsub 9(\kappa_{2.5}^{13})$.  We saw
in the preceding paragraph that $\jsub 8(\sigma_2) >
\kappa_{2.5}^{13}$, so the analysis of $\jsub 9$ done earlier implies
that $\jsub 9(\kappa_{2.5}^{13}) < \kappa_4$.  Putting all of this
together, we get
$$\align k''(k''(k)(k''(k)(k(j)(\mu)))) &=
\jsub 9\jsub {14}(\jsub 9\jsub {15}(\jsub 9\jsub {15}(\kappa_3^{10}))) \\
&< \jsub 9\jsub {14}(\jsub 9\jsub {15}(\jsub 9\jsub {15}(\kappa_{2.5}^9))) \\
&< \jsub 9\jsub {14}(\jsub 9(\kappa_{2.5}^{14})) \\
&< \jsub 9(\kappa_{2.5}^{13}) \\
&< \kappa_4, \endalign$$
as desired.

The comparisons between critical points in this section are rather
haphazard; it would be nice to have a definite algorithm for deciding
which of two given critical points is smaller.  Such an algorithm is
known~\cite{6}, but it proceeds through all smaller critical
points, and hence is useless if one wants to know whether, say,
$\jsub {12}(\jsub {13}(\kappa_2^{13}))$ lies below $\kappa_{2.5}^{10}$.
A definite question one can ask is whether there is a primitive
recursive algorithm for comparing the critical points of two
embeddings (given as expressions in $j$).

The ordinals $e({<}\alpha)$ also yield some interesting questions.
For instance, the fact that there are only finitely many modified
restrictions of members of $\jalgp$ to a given critical point~\cite{6}
means that there can only be finitely many such ordinals between two
successive critical points; in fact, there are at most $n2^n$ of them
below critical point number $n$ if the critical points are numbered in
increasing order starting with $n=0$.  But this bound is not sharp; it
turns out that there are no such ordinals between $\kappa_0$ and
$\kappa_1$, one ($\sigma_1$) between $\kappa_1$ and $\kappa_2$, one
($\jsub 3({<}\kappa_1)$) between $\kappa_2$ and~$\kappa_{2.5}$, and
five ($\jsub 7({<}\kappa_1)$, $\jsub 3(\sigma_1)$, $\jsub 2(\sigma_1)$,
$j(\sigma_1)$, and $\sigma_2$) between $\kappa_{2.5}$ and $\kappa_3$.
One can ask for improved bounds on the number of these ordinals.

A more interesting question
is whether the place at which an embedding $e$ skips over a critical
point $\beta$ (i.e., the least $\alpha$ such that $e(\alpha) \ge \beta$)
must be of the form $e'({<}\gamma)$ for some embedding $e'$ and
critical point $\gamma$.  The converse turns out to be true: if $e$
is an embedding with critical point $\gamma$, then $e'({<}\gamma)$ is
the least $\alpha$ such that $e(\alpha) \ge ee'(\gamma)$.

\head 6. Conclusion \endhead

It seems likely that the results in this paper are not optimal.  However,
pushing the lower bound on the growth rate of the number $F(n)$ of critical
points below $\kappa_n$ to a function beyond $F_{\omega+1}$ will probably
require a new idea, just as it took an additional idea to get from
$F_\omega$ to $F_{\omega+1}$.  One might hope that the growth rate can
be shown to be so great that the function $F$ cannot be proven to be
well-defined within Peano arithmetic, or perhaps even a stronger theory.
In fact, one might hope to obtain the consistency of large cardinals
from the assumption that $F(n)$ is finite for all $n$, but this would
seem to be a great deal to ask for.

The assumptions on which the construction in section~4 was based appear to
be highly coincidental, but this coincidence merely allowed us to get
a couple of levels higher in the $F_n$ hierarchy than we would have
without it.  Further experimentation with these ordinals would quite
possibly reveal more such relations, allowing for an improved
construction giving even more critical points below $\kappa_4$.

The best that one could hope for in this direction would be a version
of Figure~1 in which every critical point from the algebra $\jalgp$
would be produced by applying Lemma~2 to column~1.  As a starting
point, Laver~\cite{6} has shown that the ordinals $\gamma_m =
\crit{\jsub {2^m}}$ are in fact all of the critical points from
$\jalgp$, and that $\jsub {2^m}(\gamma_m) = \gamma_{m+1}$.  Results I
have obtained recently (a number of which were obtained independently by
A.~Dr\'apal) show that the embeddings $\jsub {2^{2^m}-1}$
have the desired form for column~$1$ in Figure~1:  \threeseq\jsub {2^{2^m}-1}
/\gamma_0/\gamma_{2^m}/\gamma_{2^{m+1}}/
for all $m$.  There is also a natural candidate for
column~$2$: \threeseq\jsub{2^{2^{2^m}}-2}/ \gamma_1/ \gamma_{2^{2^m}}/
\gamma_{2^{2^{m+1}}}/.  These calculations also show that the
first $256$ critical points from $\jalgp$ are those produced in
section~2 below $\kappa_4$, and that the next critical point beyond
these is the ordinal $\mu = \jsub 7(\kappa_2)$ used in sections~4
and~5.  This leads to some hope that the ultimate construction is
indeed possible.  In order to produce a growth rate beyond $F_\omega$,
though, this construction would probably require an additional
recursive level to produce the starting entries in whichever columns
would need them, so one can expect it to be very complex.

\enddocument